\documentclass[11pt]{article}
\usepackage{epsfig,amsmath,latexsym,amssymb}
\usepackage{amscd}
\usepackage{multirow}
\usepackage{graphicx, subcaption, caption, listings, url}
\usepackage{enumitem, color}
\usepackage[top=2.8cm,bottom=2.3cm,textwidth=16cm]{geometry}
\def\R{{\mathbb R}}  
\def\N{{\mathbb N}}  
\def\p{{\mathbb P}}  
\def\Z{{\mathbb Z}}  
\def\E{{\mathbb E}}  

\newcommand{\Remm}[1]{}
\newtheorem{theo}{Theorem}[section]
\newtheorem{lemma}[theo]{Lemma}
\newtheorem{prop}[theo]{Proposition}

\newtheorem{cor}[theo]{Corollary}
\newtheorem{defi}[theo]{Definition}
\newtheorem{model ass}[theo]{Model Assumptions}

\numberwithin{equation}{section}

\newenvironment{Proof}[1][\unskip]{\footnotesize \textbf{Proof #1.}}{ \hfill  $\square$ \par}

\begin{document}
\author{Philippe Deprez\footnote{RiskLab, Department of Mathematics, 
ETH Zurich, 8092 Zurich, Switzerland} \qquad 
Mario V.~W\"uthrich$^\ast$\footnote{Swiss Finance Institute SFI Professor, 8006 Zurich, Switzerland}}

\date{\today}
\title{
Scale-Free Percolation in Continuum Space
}
\maketitle

\begin{abstract} 
The study of real-life network modeling has become very popular in recent years. 
An attractive model is the scale-free percolation 
model on the lattice $\Z^d$, $d\ge1$, because it 
fulfills several stylized facts observed in large real-life networks. 
We adopt this model to continuum space 
which leads to a heterogeneous random-connection model on $\R^d$:
particles are generated by a 
homogeneous marked Poisson point process on $\R^d$, and  
the probability of an edge between two particles is 
determined by their marks and their distance. 
In this model we study several properties such as the degree distributions, 
percolation properties and  graph distances. 
\end{abstract}


%


\section{Introduction}
The study of real-life networks such as virtual social networks or financial networks 
has become very popular in recent years, see for example \cite{NSW, cont1, Cont}. 
Such networks can be seen as  sets of particles 
that are possibly linked to each other.  
Several stylized facts of  large real-life networks 
have been observed using large empirical data sets  
(see \cite{NSW} and Section 1.3 in \cite{Durrett} 
for further details):
\begin{itemize}
\item
The minimal number of links that connect 
two particles, called the graph distance, 
is typically small for distant particles. 
This is called the ``small-world effect''. 
There is the observation that most particles in many 
real-life networks are connected   
by at most six links, see  \cite{Watts}.
\item 
Particles that are linked tend to have 
common friends, which is called the 
``clustering property''.
\item 
The number of links of a given particle, called the degree, has a heavy-tailed 
distribution with 
(power law) tail parameter $\tau>0$. 
The tail parameter is often observed to be between $1$ and $2$, 
i.e.~the degree distribution has finite mean and infinite 
variance. We refer to \cite{Durrett} for explicit examples. 
\end{itemize}
Since it is too complicated to model large real-life 
networks particle by particle, many theoretical 
random graph models have been developed and 
their geometrical properties studied.  
One of these models is the homogeneous long-range percolation 
model on $\Z^d$, $d\ge1$, first introduced in 
\cite{Zhang} for $d=1$. 
The set of particles is the lattice $\Z^d$ and any 
two particles $x, y\in\Z^d$ are independently linked 
with probability $p_{xy}$ which behaves  as
$\lambda |x-y|^{-\alpha}$ for $|x-y|\to \infty$, with fixed 
constants $\lambda,\alpha>0$. 
Since close particles are likely linked, this   
model has a local clustering property. 
Moreover, depending on $\alpha$, the graph 
distance of two connected particles 
is roughly of logarithmic order as their separation 
tends to infinity, see \cite{biskup}. This is a 
version of the small-world effect.   
However, this model does not fulfill the stylized fact of having heavy-tailed 
degree distributions. Therefore, 
\cite{Remcoscale} extended the homogeneous long-range percolation 
model to a scale-free percolation model on $\Z^d$
(also known as inhomogeneous long-range percolation model). 
In their model they consider a collection $(W_x)_{x\in\Z^d}$ 
of i.i.d.~positive weights that are heavy-tailed with 
tail parameter $\beta>0$, and they assign to each 
particle $x\in\Z^d$ the random weight $W_x$. 
%
%
Given these weights, any two particles $x, y\in\Z^d$ are independently linked with 
probability $p_{xy}$ which is approximately 
$\lambda W_x W_y |x-y|^{-\alpha}$ for large $|x-y|$ 
and given constants $\lambda,\alpha>0$. 
Note that $p_{xy}$ is increasing in the weights $W_x$ and $W_y$, and 
decreasing in the distance between $x$ and $y$. 
This means that the weights make particles more or less attractive, 
i.e.~particles with large weights play the role of  hubs
in this network. 
This extension of the homogeneous model is 
very natural since the existence 
of hubs is often observed in real-life networks. 
Again, this model has a local clustering property. 
Depending on $\alpha$ and $\beta$, \cite{Remcoscale} showed that the 
degree distribution  is heavy-tailed, 
i.e.~this model fulfills the stylized fact of having  
heavy-tailed degree distributions. 
Moreover, they showed that whenever 
the degree distribution has finite mean 
but infinite variance, the graph distance of 
two particles behaves doubly logarithmically as their separation 
tends to infinity. This is again a version of the small-world effect, sometimes 
called the ultra-small-world effect.

~

In this article we adopt the scale-free  
percolation model on $\Z^d$ to the continuum space $\R^d$ 
as proposed in \cite{Remcoscale}, which leads 
to a heterogeneous random-connection model (RCM) on $\R^d$ 
where particles are no longer restricted to a lattice.  
Instead of taking the particles to be the vertices of $\Z^d$ 
with assigned weights, we distribute particles randomly in space according 
to a homogeneous Poisson point process on $\R^d$, and to each 
particle $x$ we attach (independently of its location) a positive random weight $W_x$ whose 
distribution is heavy-tailed with tail parameter $\beta>0$. 
Given the Poisson cloud and the weights, two 
particles $x$ and $y$ are  linked with probability $p_{xy}(\lambda,\alpha)$ 
as in the scale-free percolation model on $\Z^d$. 
This heterogeneous RCM can  be seen as an extension of the 
homogeneous RCM on $\R^d$, which was 
introduced and studied in \cite{Penrose}, 
while an applied version already appeared in \cite{Gilbert2}. 
The main reference for the homogeneous RCM and other 
continuum percolation models is~\cite{Meester-Roy}. 
The goal of this article is to prove similar results as in \cite{Remcoscale, Rajat} 
for the heterogeneous RCM. 
In particular, depending on $\alpha$ and $\beta$, we show that 
in our heterogeneous RCM the degree distribution 
is heavy-tailed with tail parameter $\tau(\alpha,\beta)>0$. 
Assuming that the weights follow a Pareto distribution with 
tail parameter $\beta$, we give an explicit 
expression of the degree distribution as well as the expected degree of a given particle 
in terms of the model parameters $\beta$, $\alpha$, $\lambda$ and 
the intensity of the Poisson point process.  
This result improves the bounds given in Proposition 2.3 in \cite{Remcoscale} 
for this particular choice of weight distribution. 
This explicit expression is also helpful to calibrate the model parameters 
for real-life network applications. 
Moreover, we show that there is a 
non-trivial phase transition depending on $\alpha$ and $\beta$, 
where $\lambda$ plays the role of the percolation parameter. 
Above criticality there is a unique infinite connected 
component. 
For real-life network applications the interesting case is $\tau(\alpha,\beta)\in (1,2)$ 
and we observe that in this case the model percolates for all $\lambda>0$.
In other words, in this latter case  
the network contains infinitely many particles 
that are all connected through links. 
We furthermore study graph distances 
between particles that lie in the same  
connected component. Similar to 
\cite{Remcoscale, Rajat} we prove 
the existence of   
different asymptotic regimes, which are characterized  
by $\alpha$ and the power law constant $\beta$ of the marks.  
As key step in that proof, we show that the size of the largest connected component 
restricted to a finite box is of the same order as the total number of particles in that box. 
This result is of independent interest and states that 
the number of particles belonging to the largest connected 
network in a finite box $[0,m)^d$ is of order $m^d$.  
%
Moreover, we show that there is no percolation at criticality 
whenever $p_{xy}(\lambda,\alpha)$ does not decrease too fast  
in the distance of two particles. 

~

Compared to inhomogeneous long-range percolation on $\Z^d$, 
the heterogeneous RCM has the advantage that  
some proofs of the results are more easy to handle  
since we can use standard integration in $\R^d$. 
This also allows us for calculating several graph properties 
explicitly which is of central interest for calibrating 
model parameters.  On the other hand, some proofs are more involved 
because one needs to make sure that the Poisson cloud 
is sufficiently regular. 
We also mention that the 
continuum space model, as an extension of the 
lattice model, has the advantage that it can  
be extended to Poisson point processes with space 
dependent (random) intensity functions. This can be used to model 
networks that have more densely populated areas than other areas.

~

The paper is organized as follows. In the next section 
we introduce the model. In Section \ref{results} we state the  
main results on the degree distributions, the percolation 
properties, the absence of percolation at criticality, 
the size of  largest connected components 
in finite boxes and the graph distances in the random graph.
Section \ref{section:degree} 
gives the proofs of the results on the degree distributions 
and in Section \ref{section:value} we prove the 
percolation properties. 
In Section \ref{Subsection: boxes} we prove 
the absence of percolation at criticality and the results 
on the size of  largest connected components 
in finite boxes are given. Finally, Section \ref{Section: Distances} 
contains the proofs of the results on graph distances.


\section{The model}
We introduce a heterogeneous RCM 
which modifies the homogeneous RCM defined
in \cite{Penrose} and which is a continuum space 
analogue to the inhomogeneous long-range percolation 
model presented by \cite{Remcoscale}. 
%
The tuple $(X,\nu, \beta,\lambda, \alpha)$  denotes a
heterogeneous RCM on $\R^d$, $d\ge1$, where we make the following
assumptions:
\begin{enumerate}
\item 
$(X,\nu,\beta)$ is a homogeneous marked Poisson point process, 
where $X$ denotes the spatially homogeneous
Poisson point process on $\R^d$ with fixed intensity $\nu>0$, 
and $\beta>0$ denotes the power law tail parameter of the 
distribution of the i.i.d.~marks $W_x$, $x\in X$. 
We assume that $W_x$, $x\in X$, has Pareto distribution 
with scale parameter $1$, i.e.  
\begin{equation*}
\p[W_x>w]=w^{-\beta}, \qquad \text{ for $w\ge 1$.}
\end{equation*}
\item 
Given $X$ and $(W_x)_{x\in X}$, we have an edge (link) between 
two distinct particles $x\neq y\in X$, write $x\Leftrightarrow y$, 
independently of all other possible edges, 
with probability
\begin{equation*}
p_{xy} =p_{xy}(\lambda, \alpha)= 1- \exp \left\{ -\lambda W_{x}W_{y}
|x-y|^{-\alpha}\right\},
\end{equation*}
with  constants $\lambda>0$ and $\alpha>0$, and $|\cdot|$ denoting 
the Euclidean norm on $\R^d$. 
\end{enumerate}
By replacing $\lambda$ by $\theta^2 \lambda$ we can extend the results
to Pareto distributions with arbitrary scale parameter $\theta>0$, 
but we use $\theta=1$ as normalization of the model. In \cite{Remcoscale} 
there is a more general version for the choice of the distribution of the marks 
$(W_x)_x$, 
but since eventually only the choice of $\beta >0$ of regular 
variation at infinity is relevant, the Pareto distribution  provides the full flavor 
of the asymptotic results. 
Moreover, for our results only the tail behavior of 
$p_{xy}$ is relevant, which is of order 
$\lambda W_x W_y |x-y|^{-\alpha}$; but we make the particular choice 
of $p_{xy}$ to simplify calculations. 
We call $X$ the Poisson cloud with particles $x\in X$. The marks
$(W_x)_{x\in X}$ are the weights in the particles $x\in X$
that determine the edge probabilities $p_{xy}$ between the corresponding particles $x$
and $y$ of  $X$.
It follows from \cite{Brug} that the model is shift invariant and ergodic.


\section{Main results}\label{results}
\subsection{Degree distribution}
We define the degree $D_x$ of particle $x\in X$ to be the number
of particles $y\in X$ such that $x$ and $y$ are linked, i.e.~$x\Leftrightarrow y$. 
Observe that the distribution of $D_x$ is translation invariant
in the sense that we may start at every particle $x$ of the Poisson cloud $X$. 
Since $D_0$ is only defined if the origin belongs to the Poisson cloud $X$,  
we consider $D_0$ under the conditional probability $\p_0$,
conditionally given that the Poisson cloud has a particle 
at the origin.  
The probability $\p_0$ is the Palm measure of $\p$, and the conditioning 
on the event of having a particle at the origin does not influence 
the rest of the Poisson process, see for instance Chapter 12 in \cite{Daley}. 
The first result describes the distribution of the degree $D_0$ 
under $\p_0$.

\begin{theo}\label{Theorem: degree}
We obtain the following cases. 
\begin{enumerate}
\item[(i)]
For $\min\{ \alpha, \beta \alpha \} \le d$ 
we obtain  $\p_0 \left[D_0=\infty \right] =1$.
\item[(ii)]
For $\min\{ \alpha, \beta \alpha \} > d$ we obtain that 
$D_0$ has (under $\p_0$) a mixed Poisson distribution with mixing distribution being 
the Pareto distribution with  
shape parameter $\tau=\beta\alpha/d>1$ and scale parameter $c_1^{1/\tau}$, where  
$c_1=c_1(d, \beta, \alpha, \lambda, \nu)
=\left(\nu v_d ~\Gamma(1-d/\alpha)
\frac{\tau }
{\tau-1}\right)^{\tau}\lambda^\beta$,
and where $v_d$ denotes the volume of the unit ball in $\R^d$. 
That is, for $k\ge0$, 
\begin{equation*}
\p_0[D_0=k] 
~=~ 
\frac{\tau c_1}{k!} \int_{c_1^{1/\tau}}^\infty t^{k-\tau-1}e^{-t}dt.  
\end{equation*}
Moreover, the survival probability of this distribution fulfills
\begin{equation*}
\lim_{n\to\infty}\frac{\p_0 \left[ D_0 > n \right]}{n^{-\tau}}
~=~c_1, 
\end{equation*}
and, hence, the degree distribution is heavy-tailed with 
tail parameter $\tau=\beta\alpha/d>1$. 
The first moment of this distribution is given by 
\begin{equation*}
\E_0 \left[ D_0\right]
=\nu v_d ~\Gamma(1-d/\alpha) \left(\frac{\tau}{\tau-1}\right)^2
\lambda^{d/\alpha}.
\end{equation*}
\end{enumerate}
\end{theo}

This theorem is the continuum space analogue to 
Theorems 2.1 and 2.2 in \cite{Remcoscale};  
the explicit expression of $\E_0 \left[ D_0\right]$ improves 
the bounds given in Proposition 2.3 in \cite{Remcoscale} 
for our choice of the distribution of the weights $(W_x)_{x \in X}$. 
We observe that if $\beta \alpha/d \le 1$ or if the decay
$|x|^{-\alpha}$ is too slow for $|x| \rightarrow \infty$, namely if $\alpha \le d$, then 
any given particle shares edges with infinitely many other particles, a.s. 
This trivial case is, of course, not of interest for real-life network modeling. 
In the non-trivial case $\min\{ \alpha, \beta \alpha \} > d$ 
the distribution of the degree of a given particle 
is heavy-tailed with tail parameter $\tau=\beta \alpha/d > 1$. 
Hence, in this latter case, the continuum space model  fulfills the stylized 
fact of having heavy-tailed degree distributions. 
This differs from the homogeneous RCM, where 
the degree distribution always is light-tailed, see formula (6.1) in \cite{Meester-Roy}. 
According to the stylized facts the interesting case for real-life applications is  
 $\tau=\beta \alpha/d \in(1,2)$ with $\alpha>d$, 
see also Section 1.4 in \cite{Durrett}.  
Note that, even if $\alpha>d$, weight distributions  
having an infinite variance ($\beta<2$) do not immediately 
imply  degree distributions having an infinite variance ($\tau <2$). 
On the other hand, under the assumption $\alpha> d$, 
if the weight distributions have a finite variance ($\beta>2$), 
the degree distributions have a finite variance ($\tau>2$) as well.

~
%
%
%
%
%
%
%

\subsection{Phase transition}
In order to study the percolation properties 
of the heterogeneous RCM, denote the  
(maximal) connected component of $x \in X$ by 
\begin{equation*}
\mathcal{C}(x)~=~\{y \in X \mid \text{there is a finite path of 
edges connecting $x$ and $y$}\},
\end{equation*}
which is the set of all particles that can be reached from $x$ 
within the network. 
The percolation probability is defined by
\begin{equation*}
\theta(\lambda)~=~\p_0\left [|\mathcal{C}(0)|=\infty \right ],
\end{equation*}
where $|\mathcal{C}(0)|$ denotes the number of particles 
in the connected component of the origin. 
The critical percolation value is defined by
\begin{equation*}
\lambda_c~=~\inf\{\lambda >0 \mid \theta(\lambda)>0\}.
\end{equation*}
By ergodicity it follows that there are only finite 
connected components, a.s., whenever $\lambda<\lambda_c$;  
and there exists an infinite connected component, a.s.,  
if $\lambda>\lambda_c$. 
By the uniqueness theorem for the homogeneous RCM, 
see Theorem 6.3 of \cite{Meester-Roy}, and the fact that $p_{xy} \in (0,1)$  
for all particles $x$ and $y$, a.s., 
such an infinite connected component is unique, a.s.;  
we denote it by $\mathcal{C}_\infty$. 

\begin{figure}
\centering
\begin{subfigure}{.45\textwidth}
\centering
\includegraphics[width=1\textwidth]{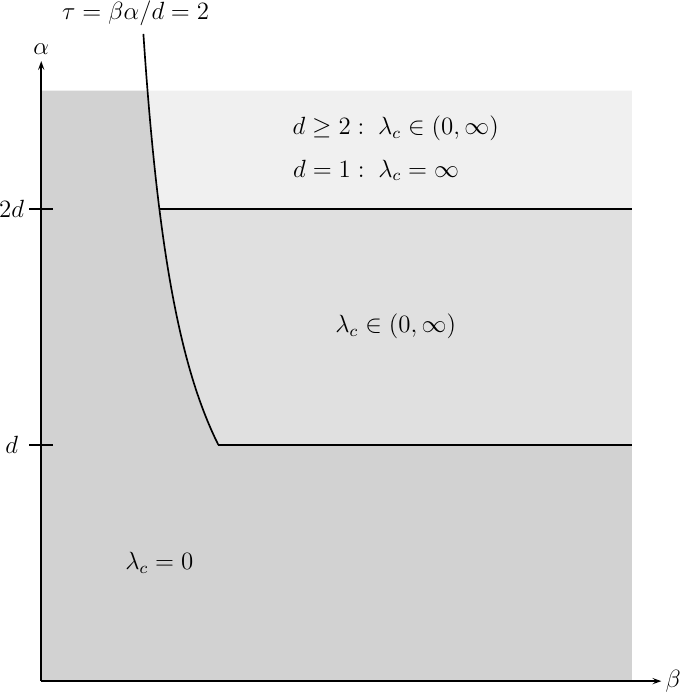}
\end{subfigure}%
\begin{subfigure}{.45\textwidth}
\centering
\includegraphics[width=1\textwidth]{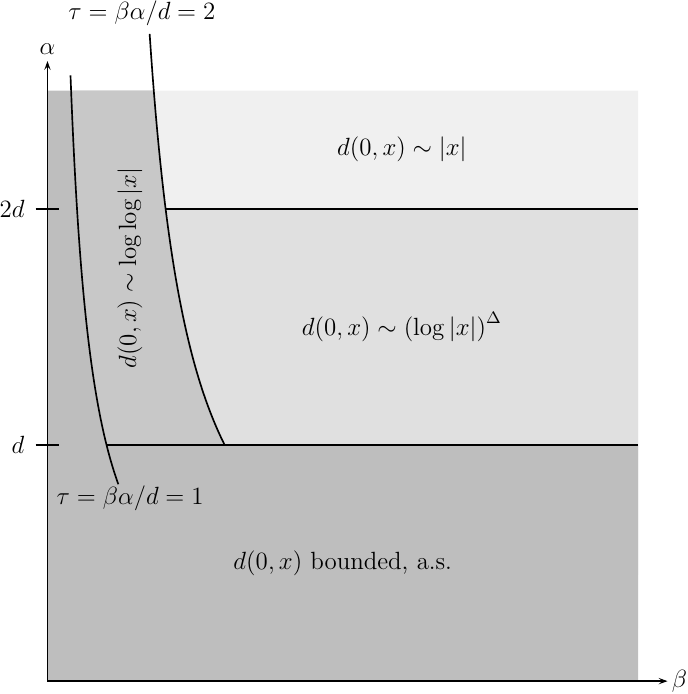}
\end{subfigure}
\caption{
The phase transition picture (lhs) and 
the graph distances (rhs) for the 
heterogeneous RCM on $\R^d$, $d\ge1$, for different model parameters $\alpha$ and $\beta$. 
Recall that the degree distributions have infinite variances if $\tau=\beta\alpha/d>2$ with $\alpha>d$. 
The existence of $\Delta$ for the asymptotic behavior of 
the graph distances in the case $\alpha\in(d,2d)$ and $\tau=\beta\alpha/d>2$ 
is still unknown; by now it is only known that $\log d(0,x) \sim \log\log |x|$. 
Additionally, in the case $\min\{\alpha,\beta\alpha\}>2d$, only the asymptotic 
linear lower bound on $d(0,x)$ is known. 
} 
\label{Picture: Phase}
\end{figure}

We refer to \cite{Bollobas1, grimmett} 
for a general introduction to percolation theory. 
For $\min\{\alpha, \beta \alpha\} \le d$ it follows from 
Theorem \ref{Theorem: degree} $(i)$  
 that $\theta(\lambda)=1$ for all $\lambda >0$,
hence $\lambda_c=0$.
The next theorem gives the percolation properties 
in the non-trivial case  $\min\{\alpha, \beta \alpha\} > d$, 
see also Figure \ref{Picture: Phase}.

\begin{theo}\label{Theorem:lambda}
Assume  $\min\{\alpha,\beta \alpha\}> d$.
\begin{enumerate}
\renewcommand{\labelenumi}{(\alph{enumi})}
\renewcommand{\labelenumii}{(\alph{enumi}\arabic{enumii})}
\item In the case  $d\ge 2$ we obtain:
	\begin{enumerate}[ref=(a\arabic*)]
	\item if $ \beta \alpha < 2d$, then $\lambda_c =0$; \label{dd:zero}
	\item if $\beta \alpha > 2d$, then $\lambda_c \in (0, \infty)$. \label{dd:zerofinite}
	\end{enumerate}
\item In the case  $d=1$ we obtain:
	\begin{enumerate}[ref=(b\arabic*)]
	\item if $ \beta \alpha < 2$, then $\lambda_c =0$;  \label{d:zero}
	\item if $\beta \alpha > 2$ and $\alpha \in (1,2]$, then $\lambda_c \in (0, \infty)$; \label{d:zerofinite}
	\item if $\min\{\alpha, \beta \alpha\} > 2$, then $\lambda_c=\infty$. \label{d:infinite}
	\end{enumerate}
\end{enumerate}
\end{theo}

This result also holds true in the discrete space model, see~\cite{Remcoscale}. 
It shows the existence of a non-trivial phase transition 
if  the degree distribution has finite variance ($\tau=\beta\alpha/d>2$) 
and $\alpha>d$ ($\alpha \in(1,2]$ in $d=1$). 
Note that in the interesting case for real-life network applications  
($\tau=\beta\alpha/d\in(1,2)$ with $\alpha>d$) there 
is a unique infinite connected component $\mathcal{C}_\infty$ for all $\lambda>0$. 
In the one-dimensional case, similar results 
to Theorem \ref{Theorem:lambda} hold true for the 
homogeneous long-range percolation model on $\Z$, where the 
probability of an edge between two sites $x$ and $y$ is given by 
$1-\exp\{-\lambda |x-y|^{-\alpha}\}$, see \cite{newman:schulman, schulman}. 
It is shown that for $\alpha \le 1$ percolation occurs for any $\lambda>0$, 
for $\alpha \in(1,2]$ percolation occurs only for $\lambda$ sufficiently large, and 
for $\alpha>2$ there does not exist an infinite connected component, a.s.

Note that $\lambda_c=0$ whenever 
$\tau=\beta\alpha/d\in(1,2)$ and, therefore, there is trivially no 
infinite connected component at criticality $\lambda_c$.  
The next theorem states that there is no infinite 
connected component at 
criticality $\lambda_c>0$ also in the case 
$\alpha\in(d,2d)$ and $\tau=\beta\alpha/d > 2$. 
This corresponds to Theorem 1.5 of \cite{Berger} 
for homogeneous long-range percolation and to 
Corollary 4 of \cite{Rajat} for inhomogeneous long-range 
percolation on the lattice. 
The case $\alpha>2d$ and $\tau=\beta\alpha/d > 2$ is still 
open, except if $d=1$ where there is never an infinite 
connected component.

\begin{theo}\label{Theorem: criticality}
Assume $\alpha \in (d,2d)$ and $\tau=\beta \alpha/d > 2$. 
There is no infinite connected component at 
criticality $\lambda_c>0$, a.s. 
\end{theo}

\subsection{Percolation on finite boxes}\label{Results:Perc}
For $n\in (0,\infty)$ we define the box $\Lambda_n=[-n,n)^d$ and 
we denote by $\mathcal{C}_n$ the 
largest connected component in $\Lambda_n$ 
(with a deterministic rule if there is more than 
one largest connected component). 
The next result shows that in case of percolation 
and $\alpha\in(d,2d)$, the number of particles 
in $\mathcal{C}_n$ presents with 
high probability at least a positive fraction 
of the Lebesgue measure of box $\Lambda_n$. 
This result is the continuum space analogue to Theorem 6 
of \cite{Rajat}.

\begin{theo}\label{Theorem:size B}
Assume $\alpha \in (d,2d)$ and $\tau=\beta \alpha/d > 1$. 
Choose $\lambda \in (0, \infty)$ such that 
$\theta (\lambda, \alpha ) > 0$. Then, for all $\alpha' \in (\alpha, 2d)$ 
there exist $\rho >0$ and $n_0 < \infty$ such that for all 
$n \ge n_0$,
\begin{equation*}
\p\left[
|\mathcal{C}_n| \ge \rho n^d
\right]
\ge
1-\exp\{-\rho n^{2d-\alpha'}\},
\end{equation*}
where $|\mathcal{C}_n|$ denotes the number of particles of the 
largest connected component in $\Lambda_n$. 
\end{theo}

For $x \in \R^d$ and $n\in(0,\infty)$ we write 
$\Lambda_n(x)=x +[-n,n)^d$ for the box of side length 
$2n$ centered at $x$. 
For $x \in X$ we write 
 $\mathcal{C}_n(x)$ for the set of 
particles in $\Lambda_n(x) \cap X$ that are connected 
to $x$ within $\Lambda_n(x)$. 
For $\ell>0$ and $\rho >0$ we call $x\in X$ a $(\rho,\ell)$-dense particle 
if the number of particles in $X$ 
belonging to $\mathcal{C}_\ell(x)$, denoted by $|\mathcal{C}_\ell(x)|$, 
is at least $\rho (2\ell)^d$, 
see also Definition 4.1 of \cite{biskup}. 
The set of $(\rho,\ell)$-dense particles in $\Lambda_n=\Lambda_n(0)$ is denoted by 
\begin{equation*}
\mathcal{D}_n^{(\rho,\ell)} 
= 
\left\{
x \in \Lambda_n\cap X\, \Big|\, |\mathcal{C}_\ell(x)| \ge \rho (2\ell)^d  
\right\}.
\end{equation*}
Corollary \ref{Corollary:3.3&3.4} below shows that whenever 
a particle $x\in X$ belongs to the infinite connected component $\mathcal{C}_\infty$, 
the probability that it is $(\rho,\ell)$-dense 
converges to $1$ as $\ell\to \infty$ for some $\rho>0$. 
This result can be interpreted as a local clustering property in the 
sense that a particle in the infinite connected component is surrounded 
by many other particles that are connected to it. 
Moreover, Corollary \ref{Corollary:3.3&3.4} shows that 
the number of dense particles in box $\Lambda_n$ 
presents with high probability at least a positive fraction of the Lebesgue measure of that box. 
Corollary \ref{Corollary:3.3&3.4} is the analogue to Corollaries 3.3 and 3.4 of \cite{biskup}. 
We use it to prove an estimate on the graph distance 
in the infinite connected component, below.

\begin{cor}\label{Corollary:3.3&3.4}
Assume $\alpha \in (d,2d)$ and $\tau=\beta \alpha/d > 1$. 
Choose $\lambda \in (0, \infty)$ such that 
$\theta (\lambda, \alpha ) > 0$.
\begin{enumerate}
\item[(i)]
There exists $\rho >0$ such that for  
$x \in \R^d$, 
\begin{equation*}
\lim_{\ell \to \infty}
\p\left[
|\mathcal{C}_\ell(x)| \ge \rho (2\ell)^d \Big | x \in \mathcal{C}_\infty
\right]
= 
1.
\end{equation*}
\item[(ii)]
For all $\alpha' \in (\alpha, 2d)$ 
there exist $\rho >0$ and $\ell_0>0$
such that for all $n>\ell_0^2$ and 
$\ell\in(\ell_0, n/\ell_0)$, 
\begin{equation*}
\p\left[
|\mathcal{D}_n^{(\rho,\ell)} | 
\ge
\rho (2n)^d
\right]
\ge
1-\exp\{-\rho n^{2d-\alpha'}\},
\end{equation*}
where $|\mathcal{D}_n^{(\rho,\ell)} |$ denotes 
the number of $(\rho,\ell)$-dense particles in $\Lambda_n$.
\end{enumerate}
\end{cor}

\subsection{Graph distances}
For $x,y \in X$ we write $d(x,y)$ for the graph distance 
or chemical distance between $x$ and $y$, i.e.~ 
\begin{equation*}
d(x,y) 
= 
\inf \left\{
n \in \N\, \big|
\exists \, x_1, \ldots , x_{n} \in X : \, 
x \Leftrightarrow x_1 \Leftrightarrow \ldots \Leftrightarrow 
x_{n-1} \Leftrightarrow x_n=y
\right\},
\end{equation*}
where we use the convention that $d(x,y)=\infty$ if 
$x$ and $y$ are not in the same connected component. 
In order to measure events involving $d(x,y)$ for 
$x,y\in\R^d$, one needs 
to make sure that the particles $x$ and $y$ lie in the 
Poisson cloud $X$. We therefore consider the $2$-fold Palm measure 
$\p_{x,y}$ of $\p$ which can be interpreted 
as the conditional distribution of the marked Poisson 
point process under the condition that there are 
particles of the process in $x$ and $y$, 
i.e.~$\p_{x,y}[\, \cdot\, ] = \p[\, \cdot\, | x,y \in X]$. 
Note that we have  
$\p[\, \cdot\, | x,y \in \mathcal{C_\infty}]=\p_{x,y}[\, \cdot\, | x,y \in \mathcal{C_\infty}]$. 
The next theorem states bounds on the graph distance 
in the case $\min\{\alpha, \beta \alpha\} >  d$, 
see also Figure \ref{Picture: Phase} 
for an illustration. 

%

\begin{theo}\label{Theorem: graph distance}
Assume $\min\{\alpha, \beta \alpha\} >  d$. 
\begin{itemize}
\item[(a)]
Assume  $\tau=\beta \alpha/d\in(1,2)$ and choose $\lambda >\lambda_c=0$. 
There exists $\eta_1>0$ such that for all $\varepsilon >0$,
\begin{equation*}
\lim_{|x| \rightarrow \infty}
\p\left[\left.
\eta_1 \frac{2}{|\log(\alpha(\beta \wedge 1)/d-1)|}
\le 
\frac{d(0,x)}{\log \log |x|}
\le 
(1+\varepsilon)
\frac{2}{|\log(\beta \alpha /d-1)|}
\right| 
0,x \in \mathcal{C}_\infty
\right]
=1.
\end{equation*}
\item[(b1)]
Assume  $\alpha \in (d,2d)$, $\tau=\beta \alpha/d>2$ and  
choose $\lambda > \lambda_c$. 
For all $\varepsilon > 0$,
\begin{equation*}
\lim_{ |x| \rightarrow \infty} 
\p\left[\left.
1-\varepsilon
\le 
\frac{\log d(0,x)}{\log \log |x|}
\le 
(1+\varepsilon) \frac{\log 2}{\log(2d/\alpha)}
\right| 0, x \in \mathcal{C}_\infty
\right]=1.
\end{equation*}
\item[(b2)]
Assume $\min\{\alpha, \beta\alpha\} > 2d$. 
There exists $\eta_2 > 0$ such that 
\begin{equation*}
\lim_{ |x| \rightarrow \infty} 
\p\left[\left.
\eta_2 
< 
\frac{d(0,x)}{|x|}
\right| 0, x \in X
\right]=1.
\end{equation*}
\end{itemize}
\end{theo}

This theorem is the continuum space analogue to 
the results in Section~5 of \cite{Remcoscale} and Theorem~8 of \cite{Rajat}. 
We  note that in homogeneous long-range percolation 
on $\Z^d$ the picture about graph 
distances in the case $\alpha>d$ (where the degree distribution has finite mean) 
has two different regimes: 
the graph distances behave roughly logarithmically if $\alpha \in (d,2d)$, 
while if $\alpha>2d$, there is a linear 
lower bound on the graph distances, 
see Theorem~1.1 of \cite{biskup} and 
Theorem~1 of \cite{Berger2}, respectively. 
In our model 
we observe the same behavior if in addition the degree 
distribution has finite variance ($\tau=\beta\alpha/d>2$). 
But if the degree distribution 
has infinite variance, we get 
an additional regime where 
the graph distances behave doubly logarithmically. 
According to the stylized facts this latter case is 
interesting for real-life network applications.  
In particular, if $d<\min\{\alpha, \beta\alpha\} < 2d$, we see 
that distant particles are connected 
by  very short paths of edges which is a version of 
the (ultra-) small-world effect. 

An upper bound on the graph distances in the case 
 $\min\{\alpha, \beta\alpha\} > 2d$ is still open. 
As in homogeneous long-range percolation with $\alpha>2d$ it is believed that 
a linear upper bound should hold, see Conjecture $1$ of \cite{Berger2}. 
For independent nearest-neighbor bond percolation on $\Z^d$, $d\ge2$, 
this result was proved in \cite{Antal}. 
Note that \textit{(b1)} states that $d(0,x)$ is roughly $(\log|x|)^\Delta$ 
for large $|x|$ and some constant $\Delta>0$. 
The existence of $\Delta$ is still unknown, even in homogeneous 
long-range percolation. 
Moreover, the optimal constants in all asymptotic behaviors are still open.

%
%
%
%
%

\section{Degree distribution}\label{section:degree}
In this section we prove Theorem \ref{Theorem: degree}. 
We start with the following observation.

\begin{lemma}\label{degree distribution}
The distribution of  degree $D_0$, conditionally given $W_0$, 
is given by
\begin{equation*}
\p_0 \left[\left. D_0=k \right| W_0\right]
=\exp \left\{-\nu \int_{\R^d}
\E_0 \left[\left.p_{0x}\right|W_0 \right] dx\right\}
\frac{\left( \nu \int_{\R^d}
\E_0 \left[\left.p_{0x}\right|W_0 \right] dx\right)^k}{k!},
\qquad \text{ for $k\in \N_0$.}
\end{equation*} 
Note that this distribution is trivial if the integral appearing twice on the right-hand side 
 does not exist. 
\end{lemma}

~

\begin{Proof}[of Lemma \ref{degree distribution}]
Let $X$ be a Poisson cloud with $0\in X$ and denote by 
$X(A)$ the number of particles in $X\cap A$ for $A\subset\R^d$. 
Conditionally given $W_0$, every particle $x\in X\setminus\{0\}$ is
now independently of the others
removed from the Poisson cloud with probability $1-p_{0x}$. 
By Proposition 1.3 of \cite{Meester-Roy}, the resulting 
process $\widetilde{X}$ is a thinned Poisson cloud, 
conditionally given $W_0$, with intensity function 
$x\mapsto \nu \E_0 \left[\left.p_{0x}\right|W_0 \right]$. 
Since $D_0=\tilde X(\R^d\setminus\{0\})$ 
in distribution, it follows that, conditionally given $W_0$,  
$D_0$ has a Poisson distribution with parameter 
$\nu \int_{\R^d}\E_0 \left[\left.p_{0x}\right|W_0 \right] dx$.
\end{Proof}

~

We now provide a 
necessary and sufficient  condition for the existence of 
$\int_{\R^d}\E_0 \left[\left.p_{0x}\right|W_0 \right] dx$ in terms of $\alpha$ and $\beta$.

\begin{prop}\label{prop:condition}
The following two statements are equivalent:
\begin{enumerate}
\renewcommand{\labelenumi}{(\roman{enumi})}
\item $\min\{\alpha, \beta\alpha\}> d$;
\item $\int_{\R^d}
\E_0 \left[\left.p_{0x}\right|W_0 \right] dx < \infty$.
\end{enumerate}
\end{prop}

Lemma \ref{degree distribution} and Proposition \ref{prop:condition} 
imply that the distribution of degree $D_0$, conditionally given $W_0$, 
has a Poisson distribution whenever $\min\{\alpha, \beta\alpha\}> d$, 
and that $D_0$ is infinite, a.s., otherwise. 
\textit{(i)} of Theorem \ref{Theorem: degree}  is therefore a 
direct consequence of Lemma \ref{degree distribution} and 
Proposition \ref{prop:condition}. 
The proof of Proposition \ref{prop:condition} is based on integral 
calculations. 

~

\begin{Proof}[of Proposition \ref{prop:condition}]
We obtain, using integration by parts in the second step, 
\begin{eqnarray*}
\E_0 \left[\left.p_{0x}\right|W_0 \right]
&=&
\int_1^\infty \beta
w^{-\beta -1}
\left(1-
\exp \left\{ -\lambda W_{0}
|x|^{-\alpha}w\right\}\right)dw
\\&=&
1-
\exp \left\{ -\lambda W_{0}
|x|^{-\alpha}\right\}
+ \lambda W_0 |x|^{-\alpha}
\int_1^\infty w^{-\beta }
\exp \left\{ -\lambda W_{0}
|x|^{-\alpha}w\right\}dw
\\&=&
1-
\exp \left\{ -\lambda W_{0}
|x|^{-\alpha}\right\}
+  \left(\lambda  W_{0} |x|^{-\alpha}\right)^{\beta }
\int_{\lambda W_0 |x|^{-\alpha}}^\infty 
z^{-\beta }
e^{ -z}dz.
\end{eqnarray*}
Note that, given $W_0$,  $1-\exp \left\{ -\lambda W_{0}|x|^{-\alpha}\right\}$ 
is integrable over $\R^d$ if and only if $\alpha>d$. It therefore 
remains to consider the integrability of 
$|x|^{-\beta\alpha}\int_{|x|^{-\alpha}}^\infty z^{-\beta }e^{-z}dz$.  
For $|x|^{\alpha}\ge1$ we obtain lower bound 
\begin{eqnarray*}
|x|^{-\beta\alpha}
\int_{|x|^{-\alpha}}^\infty 
z^{-\beta }
e^{-z}dz
~\ge~
|x|^{-\beta\alpha}
\int_1^\infty 
z^{-\beta }
e^{-z}dz,
\end{eqnarray*}
which is not integrable over $\R^d$ for $\beta \alpha \le d$.
This finally shows that {\it (ii)} implies {\it (i)}. 
For $|x|^{\alpha}\ge1$ we obtain upper bound, assume  $\beta \neq 1$, 
\begin{equation*}
|x|^{-\beta\alpha}
\int_{|x|^{-\alpha}}^\infty 
z^{-\beta }
e^{-z}dz
~\le~
|x|^{-\beta\alpha}
\left[
\int_{|x|^{-\alpha}}^1
z^{-\beta }dz
+ 
\int_1^\infty 
e^{-z}dz
\right]
~=~
|x|^{-\beta\alpha}
\left[
\frac{1-|x|^{\beta\alpha-\alpha}}{1-\beta} +e^{-1}
\right],
\end{equation*}
which is integrable over $\R^d$ for $\min\{\alpha,\beta\alpha\}>d$. 
Similarly if $\beta=1$. 
This finally shows that {\it (i)} implies {\it (ii)}. 
\end{Proof}

~

In order to prove part $(ii)$ of Theorem \ref{Theorem: degree} 
we first calculate 
$\nu\int_{\R^d}\E_0 \left[\left.p_{0x}\right|W_0 \right] dx$, 
which is finite for $\min\{ \alpha, \beta \alpha \} > d$.

\begin{prop}\label{prop 2}
Assume $\min\{ \alpha, \beta \alpha \} > d$ and set $\tau=\beta\alpha/d>1$.
We obtain
\begin{equation*}
\nu\int_{\R^d}\E_0 \left[\left.p_{0x}\right|W_0 \right] dx
~=~c_1^{1/\tau} W_{0}^{d/\alpha},
\end{equation*}
where $c_1$ is defined in Theorem \ref{Theorem: degree}, 
which has a Pareto distribution with scale parameter $c_1^{1/\tau}$ and 
shape parameter $\tau$. 
\end{prop}

\begin{Proof}[of Proposition \ref{prop 2}]
From Proposition \ref{prop:condition} we obtain that
we can apply Fubini's theorem which provides
\begin{eqnarray*}
\nu\int_{\R^d}\E_0 \left[\left.p_{0x}\right|W_0 \right] dx
&=&
\nu\int_{\R^d}\left(
\int_1^\infty \beta
w^{-\beta -1}
\left(1-
\exp \left\{ -\lambda W_{0}
|x|^{-\alpha}w\right\}\right)dw \right)dx
\\&=&
\nu\int_1^\infty
 \beta
w^{-\beta -1}
\left(\int_{\R^d}
1-\exp \left\{ -\lambda W_{0}
|x|^{-\alpha}w\right\}dx \right)dw.
\end{eqnarray*}
We first calculate the inner integral.
Using polar coordinates and integration by parts, we obtain for $w\ge 1$ 
and for $v_d$ denoting the volume of the unit ball in $\R^d$,
\begin{eqnarray*}
\int_{\R^d} 1-
\exp \left\{ -\lambda W_{0}
|x|^{-\alpha}w\right\} dx
&=& 
d v_d
\int_{0}^\infty \left(1-
\exp \left\{ -\lambda W_{0}w
r^{-\alpha}\right\}\right)r^{d-1} ~dr
\\&=&
\frac{d v_d}{\alpha}
\int_{0}^\infty \left(1-
\exp \left\{ -\lambda W_{0}w
t\right\}\right)t^{-d/\alpha -1} ~dt
\\&=&
- v_d\left(1-\exp \left\{ -\lambda W_{0}w
t\right\}\right)t^{-d/\alpha } \Big |_0^{\infty}
+ v_d \lambda W_{0}w
\int_{0}^\infty \exp \left\{ -\lambda W_{0}w
t\right\}t^{-d/\alpha} dt
\\&=&
 v_d \lambda W_{0}w
\frac{\Gamma(1-d/\alpha)}{(\lambda W_{0}w)^{1-d/\alpha}}
\int_{0}^\infty 
\frac{(\lambda W_{0}w)^{1-d/\alpha}}{\Gamma(1-d/\alpha)}
t^{1-d/\alpha-1}
\exp \left\{ -\lambda W_{0}wt\right\} dt.
\end{eqnarray*}
The latter is an integral over a gamma density for $1-d/\alpha>0$.
Therefore, we obtain
\begin{equation*}
\int_{\R^d} 1-
\exp \left\{ -\lambda W_{0}w
|x|^{-\alpha}\right\} dx
~=~
v_d ~\Gamma(1-d/\alpha) 
\left(\lambda W_{0}w\right)^{d/\alpha}.
\end{equation*}
For $\beta \alpha > d$ this implies that
\begin{equation*}
\nu\int_{\R^d}\E_0 \left[\left.p_{0x}\right|W_0 \right] dx
~=~
\nu v_d ~\Gamma(1-d/\alpha) 
\left(\lambda W_{0}\right)^{d/\alpha}\int_1^\infty
 \beta
w^{-\beta -1}
w^{d/\alpha}dw
=
\nu v_d ~\Gamma(1-d/\alpha) ~\frac{\beta}{\beta-d/\alpha}~
\left(\lambda W_{0}\right)^{d/\alpha}. 
\end{equation*}
Since $\beta/(\beta-d/\alpha)=\tau/(\tau-1)$, 
the claim follows. 
\end{Proof}

~

%

\begin{Proof}[of Theorem \ref{Theorem: degree}]
The proof of part $(i)$ is a direct consequence of 
Lemma \ref{degree distribution} and Proposition \ref{prop:condition}. 
In order to prove part $(ii)$,  
assume $\min\{\alpha,\beta\alpha\}>d$ and set $\tau=\beta\alpha/d>1$. 
Then, by Lemma \ref{degree distribution} and Proposition \ref{prop 2}, 
$D_0$ has (under $\p_0$) a mixed Poisson distribution with mixing distribution being 
the Pareto distribution with scale parameter $c_1^{1/\tau}$ and 
shape parameter $\tau$. 
From this we obtain, since 
$\E_0[W_{0}^{d/\alpha}]=\beta/(\beta-d/\alpha)=\tau/(\tau-1)$, 
\begin{equation*}
\E_0 \left[ D_0\right]
=\E_0\left[c_1^{1/\tau} W_{0}^{d/\alpha}\right]
=\nu 
v_d ~\Gamma(1-d/\alpha) \left(\frac{\tau}{\tau-1}\right)^2
\lambda^{d/\alpha}<\infty,
\end{equation*}
and for $n\ge0$, see for instance Lemma 3.1.1 of~\cite{MixedPoisson}, 
\begin{eqnarray*}
\p_0[D_0>n]
&=&
\frac{1}{n!}\int_0^\infty x^n e^{-x} \p_0\left[c_1^{1/\tau} W_{0}^{d/\alpha}>x\right]dx
\\&=&
\frac{1}{n!}\int_0^{c_1^{1/\tau}} x^n e^{-x}dx
+ 
\frac{c_1}{n!}\int_{c_1^{1/\tau}}^\infty x^{n-\tau} e^{-x} dx
~=~
e^{-c_1^{1/\tau}}\sum_{j=n+1}^\infty\frac{c_1^{j/\tau}}{j!}
+ 
\frac{c_1}{\Gamma(n+1)}\int_{c_1^{1/\tau}}^\infty x^{n-\tau} e^{-x} dx. 
\end{eqnarray*}
Choose $n\ge0$ with $n-\tau+1>0$. 
Then, 
\begin{eqnarray*}
\p_0[D_0>n]
&=&
\p_0\left[ Z >n \right] 
+ 
c_1\frac{\Gamma(n+1-\tau)}{\Gamma(n+1)} \p_0\left[  Y_{n} > c_1^{1/\tau} \right],
\end{eqnarray*}
where $Z$ has a Poisson distribution with parameter $c_1^{1/\tau}$, and 
$Y_n$ has a gamma distribution with shape parameter $n-\tau+1>0$ and scale parameter $1$.
Markov's inequality provides 
\begin{equation*}
n^\tau\p_0\left[Z\ge n \right]
~=~
n^\tau\p_0\left[e^{Z-n}\ge 1 \right]
~\le~
n^\tau\E_0\left[e^{Z}\right]e^{-n},
\end{equation*}
which converges to $0$ as $n\to\infty$. 
Moreover, by Stirling's formula, $n^\tau\Gamma(n+1-\tau)/\Gamma(n+1)$ tends to $1$ as $n\to\infty$, 
and so does $\p_0\left[  Y_{n} > s \right]$ for any $s>0$. 
Therefore, 
\begin{eqnarray*}
\lim_{n\to\infty}
n^\tau \p_0[D_0>n]
~=~
c_1.
\end{eqnarray*}

\end{Proof}

~


\section{Phase transition}\label{section:value}
In this section we prove Theorem \ref{Theorem:lambda} which 
gives the phase transition picture of the heterogeneous RCM in the 
case  $\min\{\alpha, \beta \alpha\} > d$. 
We first prove \textit{\ref{dd:zero}} and \textit{\ref{d:zero}}, namely that in any dimension 
$d\ge 1$ the critical percolation value $\lambda_c$ equals $0$  
whenever $\alpha >d$ and $\tau=\beta \alpha/d \in(1,2)$.

~

\begin{Proof}[of \textit{\ref{dd:zero}} and \textit{\ref{d:zero}}]
The idea of the proof is similar 
to the proof of Theorem 4.4 in \cite{Remcoscale}. 
Assume $\alpha >d$ and $\tau=\beta \alpha/d \in(1,2)$. 
The goal is to prove that $\theta(\lambda)>0$ for all $\lambda>0$. 
Choose $\lambda >0$ and $0<\varepsilon<\min\{d/\beta, \alpha(2/\tau-1)\}$. 
Define for $k\ge 0$ boxes $\Lambda_{2^{k}}=[-2^{k},2^{k})^d$ 
and for $k\ge 1$  define  disjoint annuli 
$R_{k}=\Lambda_{2^{k}} \setminus \Lambda_{2^{k-1}}$. 
For $k\ge1$ denote by $z_k$ the particle 
with maximal weight in $R_k$ (if it exists). 
Using that for $k\ge1$, the number of particles in $R_{k}$, denoted by $X(R_{k})$, 
has a Poisson distribution with parameter of order $\nu2^{dk}$, 
one derives that the event 
$\{X(R_{k}) \ge 1 \text{ and } W_{z_k} \ge 2^{k(d/\beta-\varepsilon)} \text{ for all } k\ge1\}$ 
has positive probability. 
Given this event, we obtain that $z_{k-1}\Leftrightarrow z_k$ 
for all $k\ge1$ with positive probability, where we set $z_0=0$. 
This implies $\theta(\lambda)>0$. We refer to \cite{Remcoscale} for the details. 
\end{Proof}

~


Next, we prove \textit{\ref{dd:zerofinite}} and \textit{\ref{d:zerofinite}}, 
which provide the non-trivial phase transition for appropriate choices 
of $\alpha$ and $\beta$. 

~

\begin{Proof}[of  \textit{\ref{dd:zerofinite}} and \textit{\ref{d:zerofinite}}]
We first show that for any dimension $d\ge 1$, 
$\lambda_c > 0$ whenever $\alpha>d$ and $\tau=\beta \alpha/d > 2$. 
We mimic the proof of Theorem 4.2 of \cite{Remcoscale}. 
Assume $\alpha>d$ and $\tau=\beta \alpha/d > 2$. 
Set $x_0=0$. 
We say that $(x_1, \ldots , x_n) \in X^{n}$ is 
a self-avoiding path in $X$ of length $n\in\N$ starting from the origin, 
write $(x_1, \ldots , x_n)$ s.a.,
if for all $i=1, \ldots , n$ 
there is an edge between $x_{i-1}$ and $x_i$, and every particle $x_i$ 
in that path occurs at most once. 
Since the degree distribution has finite mean, 
see Theorem \ref{Theorem: degree},
we obtain that the degree of each particle in such a path is bounded, a.s.
Therefore, the event that the origin lies in an infinite connected component
implies that for each $n\in\N$ there is a self-avoiding path in $X$ of 
length $n$ starting from the origin. 
Therefore, using $1-e^{-x} \le x \wedge 1$, 
\begin{equation*}
\theta(\lambda) 
\le
\E_0\left[
\sum_{\text{$(x_1, \ldots, x_{n})$ s.a.}}
\prod_{i=1}^n p_{x_{i-1}x_i}
\right]
\le
\E_0\left[
\sum_{\text{$(x_1, \ldots, x_{n})$ s.a.}}
\E_0\left[\left. 
\prod_{i=1}^{n}\left(
\frac{\lambda W_{x_{i-1}}W_{x_i}}{|x_{i-1}-x_i|^\alpha} \wedge 1\right) 
\right| X
\right]
\right].
\end{equation*}
For $n$ even and distinct $x_1,\ldots,x_{n} \in\R^d$,  
Cauchy-Schwarz' inequality implies that   
\begin{eqnarray*}
\E_0\left[
\prod_{i=1}^{n}\left(
\frac{\lambda W_{x_{i-1}}W_{x_i}}{|x_{i-1}-x_i|^\alpha} \wedge 1\right)
\right]
&=&
\E_0\left[
\prod_{i=1}^{n/2}\left(
\frac{\lambda W_{x_{2i-2}}W_{x_{2i-1}}}{|x_{2i-2}-x_{2i-1}|^\alpha} \wedge 1\right)
\prod_{i=1}^{n/2}\left(
\frac{\lambda W_{x_{2i-1}}W_{x_{2i}}}{|x_{2i-1}-x_{2i}|^\alpha} \wedge 1\right)
\right]
\\&\le&
\prod_{i=1}^{n}
\E_0\left[\left(
\frac{\lambda W_{x_{i-1}}W_{x_i}}{|x_{i-1}-x_i|^\alpha} \wedge 1
\right)^2
\right]^{1/2},
\end{eqnarray*}
and similarly for $n$ odd. 
It follows that for all $n \in \N$, 
\begin{equation*}
\theta(\lambda)
~\le~
\E_0\left[
\sum_{\text{$(x_1, \ldots, x_{n})$ s.a.}}
\prod_{i=1}^{n}
\E_0\left[\left. \left(
\frac{\lambda W_{x_{i-1}}W_{x_i}}{|x_{i-1}-x_i|^\alpha}\right)^2 \wedge 1
\right| X \right]^{1/2}
\right].
\end{equation*}
Similar to  Lemma 4.3 of \cite{Remcoscale} we get 
for $u\ge1$ and two i.i.d.~random variables $W_1$ and $W_2$  having a 
Pareto distribution with scale parameter $1$ and shape parameter $\beta$, 
using integration by parts in the first step, 
\begin{eqnarray*}
\E\left[ ( W_1W_2/u)^2 \wedge 1 \right]
&=&
\frac{1}{u^2} + \frac{2}{u^2} \int_1^u v\p\left[ W_1W_2 > v \right] dv
~=~
\frac{1}{u^2} + \frac{2}{u^2} \int_1^u v^{1-\beta}(1+\beta \log v)dv
\\&\le&
(1+\beta \log u)\left(u^{-(\beta\wedge2)} + \frac{2}{u^2} \int_1^u v^{1-\beta}dv\right)
\\&\le&
(1+1_{\{\beta\neq2\}}2/|\beta-2|) \left( 1+ \max\{2,\beta \} \log u\right )^2 u^{-(\beta \wedge 2)}, 
\end{eqnarray*}
where the last step follows by considering the cases 
$\beta < 2$, $\beta=2$ and $\beta>2$ separately. 
We finally have for any $u\ge0$, 
set $c_2=c_2(\beta)=(1+1_{\{\beta\neq2\}}2/|\beta-2|)^{1/2}$, 
\begin{equation}\label{Equation: function g}
\E\left[ ( W_1W_2/u)^2 \wedge 1 \right]^{1/2}
~\le~
1_{\left\{ u<1 \right\}} 
+
1_{\left\{ u\ge1 \right\}} 
 c_2
\left( 1+ \max\{2,\beta \} \log u\right ) u^{-(\beta/2 \wedge 1)}
~=~
g(u), 
\end{equation}
where the  equality defines the function $g$. 
Using this bound, we obtain
\begin{eqnarray}
\theta(\lambda)
&\le&\notag
\E_0\left[
\sum_{\text{$(x_1, \ldots, x_{n})$ s.a.}}
\prod_{i=1}^{n}
g(\lambda^{-1}|x_{i}-x_{i-1}|^\alpha)
\right]
\\&\le&\label{Equation: g}
\nu^{n}
\int_{\R^d}\cdots\int_{\R^d}
\prod_{i=1}^{n}
g(\lambda^{-1}|x_{i}-x_{i-1}|^\alpha)
dx_1 \cdots dx_{n}
~=~
\left(
\nu
\int_{\R^d}
g(\lambda^{-1}|x|^\alpha)dx
\right)^n.
\end{eqnarray}
Choose $\lambda<1$ such that 
$\left( 1+ \max\{2,\beta \} \log (\lambda^{-1})\right ) \lambda^{(\beta/4) \wedge (1/2)} \le 1$. 
This implies for all $u\ge0$, 
\begin{eqnarray*}
g(\lambda^{-1}u)
&=&
1_{\left\{ u<\lambda \right\}} 
+
1_{\left\{ u\ge\lambda \right\}} 
c_2
\left( 1+ \max\{2,\beta \} \log (\lambda^{-1}u)\right ) (\lambda^{-1}u)^{-(\beta/2 \wedge 1)}
\\&\le&
1_{\left\{ u<\lambda \right\}} 
+
1_{\left\{ u\ge\lambda \right\}} 
\lambda^{(\beta/4) \wedge (1/2)}
c_2
\left( 1+ \max\{2,\beta \} \log u\right ) u^{-(\beta/2 \wedge 1)}. 
\end{eqnarray*}
This provides upper bound 
\begin{equation*}
\int_{\R^d}g(\lambda^{-1}|x|^\alpha)dx
~\le~
v_d \lambda^{d/\alpha}
+
\lambda^{(\beta/4) \wedge (1/2)}
c_2
\int_{|x|\ge \lambda^{1/\alpha}}
\left( 1+ \max\{2,\beta \} \log (|x|^\alpha)\right ) |x|^{-\alpha(\beta/2 \wedge 1)}
dx.
\end{equation*}
Note that the latter integral is finite because $\alpha>d$ 
and $\beta\alpha/2 >d$. Therefore, we can choose $\lambda>0$ 
so small that the right-hand side is less than $1/(2\nu)$. 
We finally obtain for all $\lambda>0$ sufficiently small,
\begin{equation*}
\theta(\lambda)\le 2^{-n} \to 0, \qquad \text{as $n\to\infty$.}
\end{equation*}
This finally implies that for any dimension $d\ge 1$, 
$\lambda_c > 0$ whenever $\alpha>d$ and $\tau=\beta \alpha/d > 2$.

~

To finish the proof of \textit{\ref{dd:zerofinite}} and \textit{\ref{d:zerofinite}} 
it remains to show that the critical percolation value $\lambda_c$
is  finite in the case $\tau=\beta \alpha/d > 2$ and 
$\alpha >d$ ($\alpha \in(1,2]$ if $d=1$). 
We adapt the proof of Theorem 3.1 of \cite{Remcoscale}.

Partition the space $\R^d$ into cubes of side-length $n$ and 
let $r=r(n,d)$ be the maximal possible distance between two particles in neighboring cubes.
We call a cube $\Lambda$ in the partition of $\R^d$ good if $X(\Lambda)\ge 1$,
i.e.~at least one particle $x$ of the Poisson cloud $X$ falls into $\Lambda$.
Note that a cube $\Lambda$ is good with probability $1-\exp\{-\nu  n^d\}$. 
If two neighboring cubes are both good, then the probability that
the two particles with maximal weight in the respective cubes are connected is bounded below by
$1-\exp\{ -\lambda r^{-\alpha}\}$, note that $W_0\ge 1$, a.s.
We now consider the site-bond percolation model on $\Z^d$ where sites
are alive independently with probability
$1-\exp\{-\nu n^d\}$
and edges are added independently between alive nearest-neighbor 
sites with probability $1-\exp\{ -\lambda r^{-\alpha}\}$.
Note that $1-\exp\{-\nu n^d\}$ can be chosen arbitrarily
close to $1$ by taking $n$ large and then the nearest-neighbor edges can have
probabilities arbitrarily close to $1$ by choosing $\lambda$ large.
Therefore, in dimensions $d\ge 2$, the site-bond percolation model 
on $\Z^d$ percolates for sufficiently large $\lambda$, see \cite{Liggett} and 
Theorem 3.2 of \cite{Remcoscale}.
In the case $d=1$ and $\alpha \in (1,2]$ it follows from 
Theorem 1.2 of \cite{newman:schulman} that the described site-bond percolation
model percolates.
But this immediately implies that there exists an infinite connected component, a.s., 
in our model for sufficiently large $\lambda$.
\end{Proof}

~

In order to complete the proof of Theorem \ref{Theorem:lambda} it remains to 
prove \textit{\ref{d:infinite}}, namely that the critical percolation 
value $\lambda_c$ is infinite in dimension $d=1$ whenever 
$\min\{\alpha, \beta \alpha\}>2$.
The proof is similar to the one of part (c) in Theorem 3.1 of \cite{Remcoscale}. 
Since the Poisson cloud induces an additional level of complexity, 
we prove \textit{\ref{d:infinite}} in detail. 

~

\begin{Proof}[of \textit{\ref{d:infinite}}]
Assume $d=1$ and $\min\{\alpha, \beta \alpha\}>2$. 
We describe the particles of a Poisson cloud $X=(x_i)_{i\in\Z}$ containing the 
origin as follows. 
\begin{eqnarray*}
x_0&=& 0;
\\x_i&=&\inf\{x\in X | x > x_{i-1}\}, \quad \text{for $i\in\N$;}
\\x_{-i}&=&\sup\{x\in X | x < x_{-i+1}\}, \quad \text{for $i\in\N$.}
\end{eqnarray*}
Note that $x_{i}-x_{i-1}$, $i\in\Z$, 
are i.i.d.~having an exponential distribution with 
parameter $\nu$. 
The proof of  \textit{\ref{d:infinite}} is simpler in the case where the weights $(W_x)_{x \in X}$ have finite mean. We start with this case.

~

{\it Case 1.} 
Assume $\beta>1$ so that we obtain $\E_0[W_0]<\infty$.
Choose $x\in \R$ and define the event
\begin{equation*}
A_x=\{\text{no particle $y\le x$ shares an edge with any
particle $z> x$} \}.
\end{equation*}
The aim is to prove $\p_0[A_0]>0$. Stationarity and ergodicity
then imply that $A_x$ occurs for infinitely many $x\in \R$, a.s., hence
$\lambda_c=\infty$.
We define for $n\in \N$ the event
\begin{equation*}
A_0^{(n)}= 
\{ \text{
$x_{-n+k} \not\Leftrightarrow x_{k}$ for all $k=1,\ldots,n$} \}.
\end{equation*}
We get, using independence of edges,  
\begin{equation*}
\p_0[A_0]
~=~
\p_0\left[\bigcap_{n\in\N}A_0^{(n)}\right]
~=~
\E_0\left[\prod_{n\in\N}\prod_{k=1}^n
\exp\left\{-\lambda W_{x_{-n+k}}W_{x_{k}} (x_{k}-x_{-n+k})^{-\alpha}\right\}
\right].
\end{equation*}
Note that the weights $(W_x)_{x\in X}$ are independent 
for different particles and 
also independent of their locations. Using Jensen's inequality, 
we therefore get 
\begin{equation*}
\p_0[A_0]
~\ge~
\exp\left\{-\lambda
\sum_{n\in\N}\sum_{k=1}^n
\E_0\left[W_{x_{-n+k}}\right]\E_0\left[W_{x_{k}}\right]
\E_0\left[(x_{k}-x_{-n+k})^{-\alpha}\right]
\right\}.
\end{equation*}
Since for each $n\in\N$ and $k=1,\ldots, n$, 
the difference $x_{k}-x_{-n+k}$ is the sum of $n$ i.i.d.~random 
variables having  exponential distributions with parameter $\nu$, 
$x_{k}-x_{-n+k}$ has a gamma distribution with shape parameter $n$ and scale parameter $\nu$. 
Therefore,
\begin{equation*}
\p_0[A_0]
~\ge~
\exp\left\{-\lambda\E_0\left[W_0\right]^2
\sum_{n\in\N}\sum_{k=1}^n
\nu^\alpha\frac{\Gamma(n-\alpha)}{\Gamma(n)}
\right\}
~=~
\exp\left\{-\lambda\E_0\left[W_0\right]^2\nu^\alpha
\sum_{n\in\N}
n^2\frac{\Gamma(n-\alpha)}{n!}
\right\}.
\end{equation*}
By Stirling's approximation and since $\alpha>2$, the latter sum is finite, 
which finishes the proof in the case $\beta>1$. 

~

{\it Case 2.}
Assume $\beta\le1$. Using independence of edges, we obtain 
\begin{equation*}
\p_0[A_0]
~=~
\E_0\left[
\exp\left\{-\lambda 
\sum_{i,j\ge0,\, (i,j)\neq(0,0)}
W_{x_{-i}}W_{x_{j}} (x_{j}-x_{-i})^{-\alpha}\right\}
\right].
\end{equation*}
Since we condition on having a particle at the origin, 
we obtain $x_j-x_{-i}=x_j+|x_{-i}|$ 
for all $i,j\ge0$ and $(x_j+|x_{-i}|)^2\ge x_j|x_{-i}|>0$ 
for all $i,j\ge1$. This implies 
\begin{eqnarray*}
\p_0[A_0]
&=&
\E_0\left[
\exp\left\{
-\lambda W_{0}\sum_{j\ge1}W_{x_{j}} x_{j}^{-\alpha}
-\lambda W_{0}\sum_{i\ge1}W_{x_{-i}} |x_{-i}|^{-\alpha}
-\lambda 
\sum_{i,j\ge1}
W_{x_{-i}}W_{x_{j}} (x_{j}-x_{-i})^{-\alpha}\right\}
\right]
\\&\ge&
\E_0\left[
\exp\left\{
-\lambda W_{0}\sum_{j\ge1}W_{x_{j}} x_{j}^{-\alpha}
-\lambda W_{0}\sum_{i\ge1}W_{x_{-i}} |x_{-i}|^{-\alpha}
-\lambda 
\sum_{i\ge1}
W_{x_{-i}}|x_{-i}|^{-\alpha/2}
\sum_{j\ge1}
W_{x_{j}}x_{j}^{-\alpha/2}
\right\}
\right].
\end{eqnarray*}
In order to prove that all sums on the right-hand side are finite, a.s., 
it suffices to check that 
$\E_0\left[\sum_{j\ge1}\!W_{x_{j}}x_{j}^{-\alpha/2} \right]<\infty$. 
Since $\beta\le1$ and $\beta\alpha>2$ we can 
choose $\varepsilon >0$ such that 
$-\alpha/2+(1-\beta)(1-\varepsilon)/\beta<-1$,  
and we set $a_j=j^{(1+\varepsilon)/\beta}$ for $j\ge1$. 
Write 
\begin{equation}\label{Equation: two sums}
\sum_{j\ge1}W_{x_{j}}x_{j}^{-\alpha/2}
~=~
\sum_{j\ge1}\frac{W_{x_j}\wedge a_j }{x_{j}^{\alpha/2}}
+
\sum_{j\ge1}\frac{(W_{x_j}-a_j)_+}{x_{j}^{\alpha/2}}.
\end{equation}
Since 
$\p\left[ W_0 > a_j \right]= a_j^{-\beta} =j^{-(1+\varepsilon)}$, 
the Borel-Cantelli lemma implies that  $(W_{x_j}-a_j)_+$ is positive 
for only finitely many $j\ge1$, a.s. Hence, the second sum 
in \eqref{Equation: two sums} is finite, a.s.  
For the first sum in \eqref{Equation: two sums} note that 
\begin{equation*}
\E_0\left[\frac{W_{x_j}\wedge a_j }{x_{j}^{\alpha/2}}\right]
=
\E_0\left[x_{j}^{-\alpha/2}\right]\E_0\left[W_{x_j}\wedge a_j \right]
\le
\nu^{\alpha/2}\frac{\Gamma(j-\alpha/2)}{\Gamma(j)}\sum_{1\le k\le a_j}\p\left[W_0> k \right]
=
\nu^{\alpha/2}\frac{j\Gamma(j-\alpha/2)}{j!}\sum_{1\le k\le a_j}k^{-\beta}. 
\end{equation*}
For $\beta<1$ this implies for an appropriate constant $c_3>0$, 
using Stirling's approximation, 
 \begin{equation*}
\E_0\left[\frac{W_{x_j}\wedge a_j }{x_{j}^{\alpha/2}}\right]
~\le~
c_3\frac{j\Gamma(j-\alpha/2)}{j!}a_j^{1-\beta}
~=~
c_3\frac{\Gamma(j-\alpha/2)}{j!}j^{1+(1-\beta)(1+\varepsilon)/\beta}
~=~
c_3j^{-\alpha/2+(1-\beta)(1+\varepsilon)/\beta}(1+o(1)), 
\end{equation*}
as $j\to\infty$. 
By the choice of $\varepsilon>0$, the right-hand 
side is summable in $j$. The same conclusion holds true in the case $\beta=1$. 
This completes the proof of \textit{(b3)}. 
\end{Proof}

\section{Percolation on finite boxes}\label{Subsection: boxes}

In this section we prove Theorems \ref{Theorem: criticality} and \ref{Theorem:size B} 
and Corollary \ref{Corollary:3.3&3.4}. 
The key result is Lemma \ref{Lemma:crucial} below 
which corresponds to Lemma 2.3 of \cite{Berger} 
in homogeneous long-range percolation on $\Z^d$.

\begin{lemma}\label{Lemma:crucial}
Assume  $\alpha \in (d,2d)$ and $\tau=\beta \alpha/d > 1$, and  
choose $\lambda \in (0, \infty)$ with $\theta(\lambda, \alpha) >0$.  
For every $\varepsilon \in (0,1)$, $\rho > 0$ and $\alpha' <2d$ there exists 
$m' \ge 1$ such that for all $m\ge m'$,  
\begin{equation*}
\p\left[ \left| \mathcal{C}_m \right | \ge \rho m^{\alpha'/2} \right] \ge 1-\varepsilon,
\end{equation*}
where $|\mathcal{C}_m|$ denotes the number of particles of the 
largest connected component in $\Lambda_m$. 
\end{lemma}

The proof of Lemma \ref{Lemma:crucial} is based on 
a renormalization technique introduced in 
\cite{Berger}. 
We explain this renormalization in detail because the 
Poisson cloud induces an additional level of complexity. 
For integers $m \ge 1$, $k\ge 0$ and $x \in m\Z^d$ we define 
the box with corner $x$ and side length $m$ and its $k$-enlargement by
\begin{equation*}
B_m(x)= x+ [0,m)^d \quad \text{and} \quad B^{(k)}_m(x)= x+ [-k,m+k)^d,
\end{equation*}
respectively. 
We write $B_m$ and $B^{(k)}_m$ if $x=0$. 
We call a set of at least $\ell\ge1$ particles in $B_m(x) \cap X$ an $\ell${\it-semi-cluster} 
if these particles are connected within its $k$-enlargement 
$B_m^{(k)}(x)$. 
%
%

For an integer valued sequence  $(a_n)_{n\in \N_0}$ 
with $a_n > 1$, $n\ge 0$, we 
define for $n\in \N$ the cube lengths 
\begin{equation*}
m_n=a_nm_{n-1} =m_0\prod_{i=1}^n a_i =\prod_{i=0}^n a_i,
\quad \text{with $m_0=a_0$.}
\end{equation*}
For $x \in m_n\Z^d$ we call $B_{m_n}(x)$ an $n${\it-stage box}. 
Note that each $n$-stage box contains $a_n^d$ of 
$(n-1)$-stage boxes $B_{m_{n-1}}(z) \subset B_{m_n}(x)$ with 
$z \in m_{n-1}\Z^d$, which we call {\it children} of $B_{m_n}(x)$. 
In the following we recursively define  the 
aliveness of $n$-stage boxes.
%
%

\begin{defi}\label{Definition:aliveness} 
Let $(a_n)_{n\in \N_0}$ be an integer valued sequence 
with $a_n > 1$, $n\ge 0$, and define $(m_n)_{n\in \N_0}$ as above. 
Choose $k\ge0$ and 
let $(\theta_n)_{n\in \N_0}$ be a real valued sequence 
with $\theta_n \in (0,1)$ for $n \ge 0$. 
\begin{itemize}
\item
For $x\in m_0\Z^d$ we say that $0$-stage box $B_{m_0}(x)$ 
is alive if it contains a $(\theta_0a_0^d)$-semi-cluster, 
i.e.~it contains at least 
$\theta_0a_0^d$ particles that are 
connected within $B^{(k)}_{m_0}(x)$. 
\item
For $n\in\N$ and $x \in m_n\Z^d$ 
we say that $n$-stage box $B_{m_n}(x)$ is alive if the event 
$A_{n,x}=
A^{(a)}_{n,x}\cap A^{(b)}_{n,x}$ occurs, where 
\begin{itemize}
\item[$$]
$A^{(a)}_{n,x}~=~\{\text{at least $\theta_{n}a_{n}^d$ children of $B_{m_n}(x)$ are alive}\}$;
\item[$$]
$A^{(b)}_{n,x}~=~\{\text{all 
$(\prod\nolimits_{i=0}^{n-1}\theta_ia_i^d)$-semi-clusters of all alive 
children of $B_{m_n}(x)$}$ \\
\phantom{aaaaaaaaaaaaaaaaaaaaaaaaaaaaaaaaaaaaaaaaaa} 
 are connected within $B^{(k)}_{m_n}(x)\}$. 
\end{itemize}
\end{itemize}
\end{defi}

For $n\in\N_0$ define $u_n=\prod_{i=0}^{n}\theta_i$ and 
note that every alive $n$-stage box $B_{m_n}(x)$ 
contains at least 
$\prod_{i=0}^n \theta_ia_i^d=m_n^du_n$ 
particles that are connected within 
$k$-enlargement $B^{(k)}_{m_n}(x)$. 
The next lemma provides a recursive lower bound for 
$p_n=\p[A_{n,x}]$, $n\in\N$ and $x \in m_n \Z^d$.

%
%

\begin{lemma}\label{Lemma:phipsi}
Assume $\alpha \in (d,2d)$ and choose $\lambda>0$. 
Let $\xi\in(\alpha/d,2)$ and $\gamma\in(0,1)$ such that 
$18\gamma>16+\xi$. 
Choose a  real valued sequence $(\theta_n)_{n\in \N_0}$ 
with $\theta_n \in (0,1)$, $n\in\N_0$, 
and an integer valued sequence $(a_n)_{n\in \N}$ 
with $a_n > 1$, $n\in\N$. 
Assume that there exists $m_0'\ge1$ such that for 
all $a_0=m_0\ge m_0'$, all $n\in\N$ and all $s\in( 2e^{-1} \nu m_{n-1}^d,2e\nu m_{n-1}^d)$, 
\begin{equation}\label{Equation: choices}
s^\gamma < m_{n-1}^d\prod_{i=0}^{n-1}\theta_i=u_{n-1}m_{n-1}^d
\quad\text{and}\quad
s^{-\xi}<\lambda\left(\sqrt{d}m_n\right)^{-\alpha}.
\end{equation}
There exist $\varphi>0$ and $m_0''\ge m_0'$ such that for every $a_0=m_0\ge m_0''$, 
$k\ge0$, and for all $n\in\N$ and $x\in m_n\Z^d$, 
\begin{equation*}
p_n 
~=~ 
\p[A_{n,x}] 
~\ge~ 
1-
\frac{1-p_{n-1}}{1-\theta_n}
- 
a_n^{2d}\left(2e^{-2\nu m_{n-1}^d(1-2/e)}+
\left(2e^{-1}\nu m_{n-1}^d\right)^{-\varphi}\right). 
\end{equation*}
\end{lemma}

~

\begin{Proof}[of Lemma \ref{Lemma:phipsi}]
Let $n\in\N$ and $x\in m_n\Z^d$. 
We obtain 
\begin{equation}\label{bound one}
1-p_n=\p \left[A^{\rm c}_{n,x}\right]= \p \left[(A^{(a)}_{n,x}\cap A^{(b)}_{n,x})^{\rm c}\right]
\le  \p \left[(A^{(a)}_{n,x})^{\rm c} \right]
+\p\left[ (A^{(b)}_{n,x})^{\rm c}\right].
\end{equation}
For the first term in \eqref{bound one},
Markov's inequality and translation invariance provide 
\begin{eqnarray*}
\p \left[
(A^{(a)}_{n,x})^{\rm c}\right]
&=&
\p \left[\sum_{B_{m_{n-1}}(z)\subset B_{m_n}(x)}
1_{A_{n-1,z}} < \theta_{n}a_{n}^d \right]
~=~
\p \left[\sum_{B_{m_{n-1}}(z)\subset B_{m_n}(x)}
1_{A_{n-1,z}^{\rm c}} > (1-\theta_{n})a_{n}^d \right]
\\&\le&
\frac{1}{(1-\theta_{n})a_{n}^d}~
\sum_{B_{m_{n-1}}(z)\subset B_{m_n}(x)}
\p \left[A_{n-1,z}^{\rm c}\right]~=~
\frac{1}{1-\theta_{n}}~
\p \left[A_{n-1,z}^{\rm c}\right]~=~\frac{1-p_{n-1}}{1-\theta_{n}}.
\end{eqnarray*}
The second term in \eqref{bound one} is more involved due
to possible dependence in the $k$-enlargements, $k\ge0$. 
For two children $B^1$ and $B^2$ of $B_{m_n}(x)$ let $E(B^1,B^2)$ be 
the event that at least two 
$(u_{n-1}m_{n-1}^d)$-semi-clusters 
in $B^1\cup B^2$ are not connected within $B_{m_n}^{(k)}(x)$. 
We obtain 
\begin{equation}\label{Equation: estimate}
\p \left[(A^{(b)}_{n,x})^{\rm c}\right]
~\le~ 
\binom{a_n^d}{2} \sup_{(B^1,B^2)}
\p\left[
E(B^1,B^2)
\right], 
\end{equation}
where the supremum is taken over all possible choices of 
two distinct children $B^1$ and $B^2$ of $B_{m_n}(x)$. 
We fix two different  children $B^1$ and $B^2$   of $B_{m_n}(x)$. 
We obtain 
\begin{eqnarray*}
\p\left[ 
E(B^1,B^2)
\right]
&\le& 
\p\left[ 
E(B^1,B^2), \, 
2e^{-1}\nu m_{n-1}^d < X(B^1\cup B^2) < 2e\nu m_{n-1}^d
\right]
\\&&\hspace{1cm}
+~ 
\p\left[X(B^1\cup B^2)\ge 2e\nu m_{n-1}^d \right]
+
\p\left[
X(B^1\cup B^2) \le 2e^{-1} \nu m_{n-1}^d
\right],
\end{eqnarray*}
where $X(A)$ denotes the number of particles of $X$  in $A\subset\R^d$. 
A random variable $Y$ having a Poisson distribution 
with parameter $\mu>0$ satisfies, using Chernoff's bound, 
\begin{equation}\label{Equation: Chernoff}
\p\left[
Y \le e^{-1} \mu
\right]
~\le~
e^{-\mu(1-2/e)}
\quad \text{and} \quad
\p\left[
Y \ge e \mu
\right]
~\le~
e^{-\mu},
\end{equation}
see for instance (A.12) of \cite{FM}. 
Using these bounds we obtain 
\begin{equation*}
\p\left[ 
E(B^1,B^2)
\right]
~\le~ 
\p\left[ 
E(B^1,B^2), \, 
2e^{-1}\nu m_{n-1}^d < X(B^1\cup B^2) < 2e\nu m_{n-1}^d
\right]
+ 
2e^{-2\nu m_{n-1}^d(1-2/e)}.
\end{equation*}
To estimate the probability above we will condition on 
the Poisson cloud restricted to $B_{m_n}^{(k)}(x)$. 
We fix integers $s\in( 2e^{-1} \nu m_{n-1}^d,2e\nu m_{n-1}^d)$, 
$t\ge0$ and choose 
$x_1, \ldots , x_{s+t} \in B_{m_n}^{(k)}(x)$ with $x_1, \ldots , x_s \in B^1\cup B^2$. 
Assume that $\{x_1, \ldots , x_s\} = (B^1\cup B^2) \cap X$
and $\{x_{1}, \ldots , x_{t}\} = B_{m_n}^{(k)}(x) \cap X$. 
Consider the edge probabilities 
\begin{equation*}
\widetilde p_{x_ix_j} ~=~1-\exp\{-\lambda|x_i-x_j|^{-\alpha}\} 
~\le~
1-\exp\{-\lambda W_{x_i}W_{x_j}|x_i-x_j|^{-\alpha}\} 
~=~
p_{x_ix_j},
\end{equation*}
a.s., for every $i \neq j \in \{1,\ldots, s+t\}$. 
Denote by 
$\widetilde\p_X$ the probability measure of the resulting 
edge configurations restricted to $\{x_1, \ldots , x_{s+t}\}$ 
induced by $\widetilde p_{x_ix_j}$, $i \neq j \in \{1,\ldots, s+t\}$. 
Note that on $B_{m_n}^{(k)}(x) \cap X = \{x_{1}, \ldots , x_{s+t}\}$, 
$E(B^1,B^2)$ is determined by 
edges with end points in $\{x_{1}, \ldots , x_{s+t}\}$. 
We therefore get 
\begin{equation*}
\p\left[ E(B^1,B^2)\bigg | 
X\cap{\left(B^1\cup B^2\right)} = \{x_1, \ldots , x_s\},\, 
X\cap{B_{m_n}^{(k)}(x)} = \{x_{1}, \ldots , x_{s+t}\}  \right]
 ~\le~
 \widetilde\p_X[E(B^1,B^2)].
\end{equation*}
We can now argue as in the proof of Lemma 2.3 of \cite{Berger11}, 
which we briefly recall. 
By an abuse of notation, assume that $B^1\cup B^2=\{x_1, \ldots , x_s\}$ and 
$B_{m_n}^{(k)}(x) = \{x_{1}, \ldots , x_{s+t}\}$. 
Using that sequence $(a_n)_{n\in\N_0}$ satisfies \eqref{Equation: choices} 
and $B^1\cup B^2 \subset B_{m_n}(x)$, 
we obtain 
$\widetilde p_{x_ix_j}>1-\exp\{-s^{-\xi}\}=\nu_n>0$ 
for all $x_i \neq x_j \in B^1\cup B^2$, where the equality defines $\nu_n$. 
For $x_i \neq x_j \in B^1\cup B^2$ we define $\widetilde q_{x_ix_j}$ by 
$\widetilde p_{x_ix_j}=\widetilde q_{x_ix_j}+\nu_n - \nu_n\widetilde q_{x_ix_j}$. 
We now sample an edge configuration $\omega$ 
on $B_{m_n}^{(k)}(x)$ induced by 
the $\widetilde p_{x_ix_j}$'s in two steps. We sample $\omega'$ 
on $B_{m_n}^{(k)}(x)$ 
with edge probabilities $\widetilde q_{x_ix_j}$ if $x_i \neq x_j \in B^1\cup B^2$ 
and with edge probabilities  $\widetilde p_{x_ix_j}$ otherwise. 
$\omega''$ is then independently sampled on $B^1 \cup B^2$ with edge probabilities $\nu_n$. 
Note that $\omega=\omega' \vee \omega''$ in distribution. 
Let $S_1 \neq S_2 \subset B^1\cup B^2$ be two disjoint maximal 
sets in $B^1 \cup B^2$ that are each connected within $B_{m_n}^{(k)}(x)$ by  
$\omega'$-edges. 
Here, the maximality of set $S_1$ means that there is no particle in $B^1 \cup B^2$ 
which does not belong to $S_1$ but is connected to $S_1$ within $B_{m_n}^{(k)}(x)$ 
by $\omega'$-edges. 
Given $\omega'$, the probability that 
there is an $\omega$-edge between $S_1$ and $S_2$ 
is equal to the probability that there is an $\omega''$-edge 
between $S_1$ and $S_2$, given $\omega'$, which follows by maximality 
of $S_1$ and $S_2$. 
The latter probability is by definition of $\omega''$ and $\nu_n$ given by 
$1-\exp\{-|S_1||S_2|s^{-\xi}\}$, 
where $|S_i|$ denotes the number of particles of $S_i$, 
$i=1,2$. 
Given $\omega'$, denoting by $S_1,\ldots,S_l$ all disjoint maximal sets  
in $B^1 \cup B^2$ that are connected within $B_{m_n}^{(k)}(x)$ by  
$\omega'$-edges, the indices $\{1,\ldots,l\}$ form an 
inhomogeneous random graph of size $\sum_{i=1}^l|S_i|=s$ 
and parameter $\xi$, as defined in \cite{Berger11}.  
Lemma 2.5 of \cite{Berger11} shows that there 
exists $\varphi>0$ such that for all $M$ sufficiently large and all 
inhomogeneous random graphs 
of size $M$ and parameter $\xi$, the probability that such a  
graph contains more than one cluster of size at least $M^\gamma$, 
is at most $M^{-\varphi}$. This implies that there exist $\varphi>0$ and 
$m_0''\ge m_0'$ such that for all $m_0\ge m_0''$, given $\omega'$, there 
is at most one cluster of size at least 
$s^\gamma$ formed by $S_1,\ldots, S_l$ that are 
connected by $\omega$-edges, with probability at most $s^{-\varphi}$.  
Note that the existence of $\varphi$ is 
uniform in $s\in( 2e^{-1} \nu m_{n-1}^d,2e\nu m_{n-1}^d)$, 
$t\ge0$, the locations of $\{x_{1}, \ldots , x_{s+t}\}$, $\omega'$, 
$n\ge 1$, $k\ge0$ and $m_0\ge m_0''$. 
Using that sequences $(a_n)_{n\in\N_0}$ and $(\theta_n)_{n\in\N_0}$ 
satisfy \eqref{Equation: choices}, i.e. $s^\gamma < u_{n-1}m_{n-1}^d$, 
 we conclude that 
there exist $\varphi>0$  and $m_0''\ge m_0'$ such that 
for every $m_0\ge m_0''$, $k\ge0$, $n\in\N$, 
$s\in( 2e^{-1} \nu m_{n-1}^d,2e\nu m_{n-1}^d)$, $t\ge0$, 
and $x_1, \ldots , x_{s+t} \in B_{m_n}^{(k)}(x)$ with $x_1, \ldots , x_s \in B^1\cup B^2$, 
\begin{equation*}
\widetilde\p_X[E(B^1,B^2)]~\le~s^{-\varphi}.
\end{equation*}
Integrating over the particles in $B_{m_n}^{(k)}(x) \setminus (B^1\cup B^2)$ 
and $B^1\cup B^2$ 
we then get for all $m_0\ge m_0''$, 
\begin{eqnarray*}
\p\left[ 
E(B^1,B^2)
\right]
&\le&
2e^{-2\nu m_{n-1}^d(1-2/e)}
+
\sum_{s=\lfloor 2e^{-1} \nu m_{n-1}^d+1\rfloor}^{\lceil2e\nu m_{n-1}^d-1\rceil}
\p\left[ X(B^1\cup B^2)=s\right]
s^{-\varphi}
\\&\le&
2e^{-2\nu m_{n-1}^d(1-2/e)}+
\left(2e^{-1}\nu m_{n-1}^d\right)^{-\varphi}
\sum_{s=\lfloor 2e^{-1} \nu m_{n-1}^d+1\rfloor}^{\lceil2e\nu m_{n-1}^d-1\rceil}
\p\left[ X(B^1\cup B^2)=s\right]
\\&\le&
2e^{-2\nu m_{n-1}^d(1-2/e)}+
\left(2e^{-1}\nu m_{n-1}^d\right)^{-\varphi}, 
\end{eqnarray*}
which together with \eqref{Equation: estimate} implies  
for all $m_0\ge m_0''$, $k\ge0$ and $n\in\N$, 
\begin{eqnarray*}
\p \left[(A^{(b)}_{n,x})^{\rm c}\right]
&\le& 
a_n^{2d}\left(2e^{-2\nu m_{n-1}^d(1-2/e)}+
\left(2e^{-1}\nu m_{n-1}^d\right)^{-\varphi}\right).
\end{eqnarray*}

This finishes the proof of Lemma \ref{Lemma:phipsi}. 
\end{Proof}

~

Using the recursion in Lemma \ref{Lemma:phipsi} 
and sequences $(a_n)_{n\in \N}=(n^a)_{n\in \N}$ and $(\theta_n)_{n\in \N}=(n^{-b})_{n\in \N}$ 
for appropriate constants $1<a<b$, 
for explicit choices see (5) in \cite{Berger11}, 
$1-p_n$ can be bounded by a multiple (independent of $n$) 
of $1-p_0$. 
In order to prove that $p_n$ is arbitrarily close to $1$, 
it remains to show that  
$0$-stage boxes are alive with arbitrarily high probability 
for sufficiently large $a_0=m_0$ and well chosen $\theta_0\in(0,1)$. 
%

%

~

\begin{Proof}[of Lemma \ref{Lemma:crucial}]
Note that $\alpha \in (d,2d)$ and $\tau=\beta \alpha/d>1$ 
imply  $\lambda_c < \infty$, see Theorem \ref{Theorem:lambda}. 
Hence, there exists 
$\lambda \in (0, \infty)$ with $\theta=\theta(\lambda, \alpha) >0$ and 
for these parameters we have a unique infinite connected 
component $\mathcal{C}_\infty \subset X$, a.s. 
In order to prove Lemma \ref{Lemma:crucial} it suffices to check that 
for any $\varepsilon' \in (0,1)$ there exists $\theta_0\in(0,1)$ 
such that $p_0 = \p\left[ A_{0,0}\right] > 1-\varepsilon'$ for all $m_0$ 
sufficiently large.
The remainder of the proof then follows from the 
one of Lemma 2.3 of \cite{Berger11} using Lemma \ref{Lemma:phipsi} 
above with appropriate sequences $(a_n)_{n\in \N_0}$ and $(\theta_n)_{n\in \N_0}$. 
Choose $\varepsilon' \in (0,1)$ and set $\theta_0=\nu v_d 2^{-d}\theta/2$, 
where $v_d$ denotes the Lebesgue measure of the unit ball in $\R^d$. 
For $m \ge 1$ denote by $\left | \mathcal{C}_\infty \cap B_{m} \right |$ 
the number of particles in $B_m$ that belong to the infinite connected component $\mathcal{C}_\infty$. 
For $m \ge 1$ it holds that 
\begin{equation*}
\p\left[ \left | \mathcal{C}_\infty \cap B_{m} \right | \ge \theta_0m^d \right] 
~\ge~
\p\left[ \sum_{x \in X \cap  B(m/2)} 1_{\{x \in \mathcal{C}_\infty\}} \ge \theta_0m^d \right], 
\end{equation*}
where $ B(m/2)$ denotes the ball of radius $m/2$ around the origin. 
Using ergodicity and $(12.4.3)$ of \cite{Daley} it follows that
\begin{equation*}
\frac{1}{\nu v_d (m/2)^d}\sum_{x \in X \cap  B(m/2)} 1_{\{x \in \mathcal{C}_\infty\}}
~\longrightarrow~ 
\p_0\left[ 0 \in \mathcal{C}_\infty \right] 
= 
\theta,
\quad \text{ as $m \to \infty$, a.s.}
\end{equation*}
Hence, for all $m_0$ sufficiently large  
we obtain  
\begin{equation*}
\p\left[ \left | \mathcal{C}_\infty \cap B_{m_0} \right | \ge \theta_0m_0^d \right] 
~\ge~
1-\varepsilon'/2. 
\end{equation*}
Since the infinite connected component $\mathcal{C}_\infty$ is unique, a.s., 
there exists $k=k(m_0)\ge 0$ such that  $\mathcal{C}_\infty \cap B_{m_0}$ is 
connected within  $k$-enlargement $B_{m_0}^{(k)}$. 
Choose $k=k(m_0)\ge 0$ such that
\begin{equation*}
\p\left[\mathcal{C}_\infty \cap B_{m_0} \text{ is connected within } B_{m_0}^{(k)} \right]
> 1-\varepsilon'/2.
\end{equation*} 
This implies that for all $m_0$ sufficiently large, 
\begin{equation*}
p_0 
~=~ 
\p\left[ A_{0,0}\right]
~\ge~
\p\left[
\left | \mathcal{C}_\infty \cap B_{m_0} \right | \ge \theta_0m_0^d \text{ and }
\mathcal{C}_\infty \cap B_{m_0} \text{ is connected within } B_{m_0}^{(k)} 
\right]
~\ge~
1- \varepsilon'.
\end{equation*}
\end{Proof}

~

\begin{Proof}[of Theorems \ref{Theorem: criticality} and \ref{Theorem:size B}]
Note that by Lemma \ref{Lemma:crucial}, 
the number of particles of the largest connected component in $\Lambda_m$ 
is at least $\rho m^{\alpha'}$ for $\rho>0$ and $\alpha'<d$, with high probability. 
In order to prove that this number is proportional to $m^d$, with high 
probability, we apply a second renormalization based on site-bond percolation, 
as in the proof of  Theorem 6 in \cite{Rajat}, 
where the bound on the probability that a site is alive is given by 
Lemma \ref{Lemma:crucial} above. 
Theorem \ref{Theorem:size B}
then follows directly from the results in \cite{Rajat}. 
This is immediately clear because 
the bounds on attachedness of alive sites derived in Lemma 10 (a) in \cite{Rajat} 
also apply to the Poisson case, see in particular estimate (2) in \cite{Rajat}. 
In the same way, using a site-bond percolation model, Theorem \ref{Theorem: criticality} 
follows from Theorem 3  in \cite{Rajat}. 
\end{Proof}

~

\begin{Proof}[of Corollary \ref{Corollary:3.3&3.4}]
Using Theorem \ref{Theorem:size B}, the 
statements $(i)$ and $(ii)$ of Corollary \ref{Corollary:3.3&3.4}
follow from Corollaries 3.3 and 3.4 of \cite{biskup}, respectively. 
Note that Lemma 3.5 therein corresponds to Lemma 10 (a) in \cite{Rajat}.  
We refer to \cite{biskup} for the details. 
\end{Proof}


\section{Graph distances}\label{Section: Distances}
  
\subsection{Infinite variance of degree distribution, upper bound}  
In order to prove Theorem \ref{Theorem: graph distance} 
we first prove the upper bound of statement \textit{(a)} which we recall in the 
next proposition.

\begin{prop}\label{Proposition: (a)1}
Assume $\alpha>  d$ 
and $\tau=\beta \alpha/d \in (1, 2)$. 
For every $\lambda>\lambda_c = 0$ and $\varepsilon >0$ 
we obtain  
\begin{equation*}
\lim_{|x| \rightarrow \infty}
\p\left[\left.
d(0,x)
\le 
(1+\varepsilon)
\frac{2\log \log |x|}{|\log(\beta \alpha /d-1)|}
\right| 
0,x \in \mathcal{C}_\infty
\right]
=1.
\end{equation*}
\end{prop}

For the proof of this proposition  we use the following technical lemma.

\begin{lemma}\label{Lemma:distance1}
Assume  $\alpha>  d$ and $\tau=\beta \alpha/d \in (1, 2)$.  
If there exists a particle in $\Lambda_r=[-r,r)^d$, $r>0$, let $z_0 \in X \cap  \Lambda_r$ 
be the particle with maximal weight $W_{z_0}$ in $ \Lambda_r$. 
We obtain for any constant $c\ge1$, 
\begin{equation*}
\lim_{w \rightarrow \infty}
\p\left[
 |\mathcal{C}(z_0)|=\infty \Big| W_{z_0} \ge w,\, X( \Lambda_r)\ge c
\right]
=1,
\end{equation*}
where $\left| \mathcal{C}(z_0) \right |$ denotes the number of 
particles of the connected component of $z_0$, and where 
$X( \Lambda_r)$ denotes the number of particles in box $\Lambda_r$. 
\end{lemma}

Note that this result corresponds to Lemma 5.2 of \cite{Remcoscale}. 
The only difference lies in the additional condition that 
$\Lambda_r$ contains at least $c\ge1$ particles which 
ensures that  particle $z_0$ with maximal weight 
in $\Lambda_r$ exists. 

~

\begin{Proof}[of Lemma \ref{Lemma:distance1}] 
Let  $\alpha>  d$ and $\tau=\beta \alpha/d \in (1, 2)$.  
Choose $b>1$ such that $2d/\beta - b\alpha>0$ 
and choose $\varepsilon \in (\alpha/2,d/\beta)$. 
Choose $r>0$ and for $w \ge 1$ sufficiently large we define disjoint annuli 
$R_1=\Lambda_{w^{1/\alpha}}\setminus \Lambda_r$, 
$R_2=\Lambda_{w^b2^{2}}\setminus \Lambda_{w^{1/\alpha}}$ and 
$R_{k}=  \Lambda_{w^b2^{k}} \setminus  \Lambda_{w^b2^{k-1}}$ 
for $k\ge 3$. 
If $X( \Lambda_r)  \ge c$, let $z_0 \in X \cap  \Lambda_r$ 
be the particle with maximal weight $W_{z_0}$ in $\Lambda_r$.  
For $k\ge 1$, if $X(R_{k}) \ge 1$, let $z_k \in X\cap R_{k}$ be 
the particle with maximal weight $W_{z_k}$ in $R_{k}$. 
By the choices of $b$ and $\varepsilon$, and  
since $X(R_{k})$ is Poisson distributed, we obtain that the event 
$\{X(R_{k}) \ge 1 \text{ and } W_{z_k}\ge2^{k\varepsilon}w^{d/\beta} \text{ for all } k\ge 1\}$ 
has probability arbitrarily close to $1$ for $w$ sufficiently large. 
Given this event and the event 
$\{X( \Lambda_r)\ge c,\, W_{z_0} \ge w\}$, 
we obtain 
that $z_0\Leftrightarrow z_{1}$ and $z_{k+1} \Leftrightarrow z_{k}$ for all 
$k\ge1$ with probability 
arbitrarily close to $1$ for $w$ sufficiently large, 
which implies the claim. 
We refer to the corresponding proof in \cite{Remcoscale} for the details. 
\end{Proof}

~

\begin{Proof}[of Proposition \ref{Proposition: (a)1}]
Let $\alpha>  d$ 
and $\tau=\beta \alpha/d \in (1, 2)$. 
Fix $\varepsilon>0$ and choose $b\in (0,1)$ such that
\begin{equation*}
d(1+b)/\beta > \alpha
\quad \text{and}\quad
\frac{1+\varepsilon/2}{|\log b|}
\le
\frac{1+\varepsilon}{|\log(\beta\alpha/d -1)|}.
\end{equation*}
Fix $m>\max\{e,3^{1/(1-b)},c_4^{-1/d},m_0,m_1\}$, 
with $c_4=(1-2/e)\nu2^{-d}$ and 
where $m_0=m_0(b)\ge0$ and $m_1=m_1(b)\ge0$ are defined below. 
Choose $x \in \R^d$ with $|x| \ge e^{(\log m)/b}$ and set 
\begin{equation*}
k=k(x)
=
\left \lfloor \frac{\log\log|x| - \log \log m}{|\log b|} \right \rfloor \ge 1.
\end{equation*}
Note that $m\le|x|^{b^k}\le m^{1/b}$ for all $x\in\R^d$. 
For $i=0,1, \ldots, k$ write $ \Lambda(x, b^i)$ for the box 
with side length $|x|^{b^i}/2$ centered at the point at 
distance $|x|^{b^i}/2$ from the origin on the segment 
with end points $0$ and $x$. 
If there is a particle in $ \Lambda(x, b^i)$, let $z_i \in  \Lambda(x, b^i) \cap X$ 
be the particle with maximal weight $W_{z_i}$ in $ \Lambda(x, b^i)$. 
Moreover,  write $ \Lambda'(x, b^i)$ for the box 
with side length $|x|^{b^i}/2$ centered at the point at 
distance $|x|^{b^i}/2$ from $x$ (instead of the origin) 
on the segment with end points $0$ and $x$. 
 If there is a particle in $ \Lambda'(x, b^i)$, let
$z'_i \in   \Lambda'(x, b^i) \cap X$ be the particle 
with maximal weight $W_{z'_i}$ in $ \Lambda'(x, b^i)$, 
this choice is similar to the one in the proof of Theorem 5.1 of \cite{Remcoscale}. 
In order to make sure that  particles $z_i$ and $z_i'$ exist 
and that boxes  $\Lambda(x, b^i)$ and  $\Lambda'(x, b^i)$ contain 
sufficiently many particles 
for all $i=0, \ldots, k$ we consider the probability 
measure 
\begin{equation*}
\p^{k}[\, \cdot \,] 
= 
\p\left[ \, \cdot \,\Big | \, X( \Lambda(x, b^i)) \ge c_4 |x|^{db^i} 
\text{ and }  X( \Lambda'(x, b^i)) \ge c_4 |x|^{db^i} 
\text{ for all } i=0, \ldots , k
\right],
\end{equation*}
which is 
the conditional probability given that there are at least 
$c_4 |x|^{db^i}\ge c_4m^d\ge1$ particles in each of 
boxes $ \Lambda(x, b^i)$ and $ \Lambda'(x, b^i)$ for $i=0, \ldots, k$. 
Using Chernoff's bound, see \eqref{Equation: Chernoff}, 
we obtain for each $i=0, \ldots, k$, since $X( \Lambda(x, b^i))$ 
has a Poisson distribution with parameter $\nu 2^{-d} |x|^{db^i}$, 
\begin{equation*}
\p\left[
X( \Lambda(x, b^i)) < c_4 |x|^{db^i}
\right]
\le
\p\left[
X( \Lambda(x, b^i)) \le e^{-1}\nu 2^{-d} |x|^{db^i}
\right]
\le
e^{-c_4|x|^{db^i}},
\end{equation*}
for each $i=0,\ldots, k$. 
Note that $m>3^{1/(1-b)}$ and $m\le|x|^{b^k}$ imply 
$|x|^{b^i}-3|x|^{b^{i+1}}>0$ for all $i=0, \ldots, k-1$. 
This implies that all boxes $ \Lambda(x, b^i)$ are disjoint for $i=0, \ldots, k$.
It follows that
\begin{equation*}
\p\left[
X( \Lambda(x, b^i)) \ge c_4 |x|^{db^i} 
\text{ and }  X( \Lambda'(x, b^i)) \ge c_4 |x|^{db^i} 
\text{ for all } i=0, \ldots , k
\right]
\ge
\prod_{i=0}^k
\left(1-e^{-c_4|x|^{db^i}}\right)^2,
\end{equation*}
where the inequality comes from the fact that 
$\Lambda(x, 1)=\Lambda'(x, 1)$. 
We write $k=k(x)=\lfloor N(x) - M \rfloor$ with  $N=N(x)=(\log \log |x|) / |\log b|$ 
and $M=(\log \log m) / |\log b|$, and we obtain,  
note that $|\log b|=-\log b$  (because $b\in(0,1)$), 
\begin{eqnarray}
\lim_{|x| \to \infty }\sum_{i=0}^{k}
\exp\left\{
-c_4
|x|^{db^{i}}
\right\}
\!\!\!&=&\!\!\! \notag
\lim_{N \to \infty }\sum_{i=0}^{\lfloor N-M \rfloor}
\exp\left\{
-c_4
e^{b^{-N} db^{i}}
\right\}
\le
\lim_{N \to \infty }\sum_{i=0}^{\lfloor N-M \rfloor}
\exp\left\{
-c_4
e^{db^{-(\lfloor N-M \rfloor-i)-M}}
\right\}
\\&=&\!\!\! \label{Equation: similar}
\lim_{N \to \infty }\sum_{ j=0}^{\lfloor N-M \rfloor}
\exp\left\{
-c_4
e^{db^{-j-M}}
\right\}
=
\sum_{ j \ge 0}
\exp\left\{
-c_4
e^{db^{-j-M}}
\right\} \in (0,\infty).
\end{eqnarray}
For fixed $\varepsilon'>0$ we therefore can choose $m_0\ge0$ 
(used for the choice of $m>m_0$) so large that 
for any sufficiently large $|x|$, 
the first equality defines event $N_k=N_k(x,b)$, 
\begin{equation}\label{Equation: kprob}
\p\left[
N_k
\right]
=
\p\left[
X( \Lambda(x, b^i)) \ge c_4 |x|^{db^i} 
\text{ and }  X( \Lambda'(x, b^i)) \ge c_4 |x|^{db^i} 
\text{ for all } i=0, \ldots , k
\right]
\ge
1-\varepsilon'.
\end{equation}
Note that for every $\delta \in (0,1)$ and $i=0,\ldots,k$, 
\begin{eqnarray*}
\p^k\left[
W_{z_i} \le \left(c_4 |x|^{db^i}\right)^{(1-\delta)/\beta}
\right]
&\le&
\p\left[\left.
W_{z_i} \le X( \Lambda(x, b^i))^{(1-\delta)/\beta}
\right| X( \Lambda(x, b^i)) \ge c_4 |x|^{db^i}
\right]
\\&=&
\p\left[\left.
\max_{z \in  \Lambda(x, b^i) \cap X} W_z \le X( \Lambda(x, b^i))^{(1-\delta)/\beta}
\right| X( \Lambda(x, b^i)) \ge c_4 |x|^{db^i}
\right]
\\&=&
\E\left[\left.
\left(
1- X( \Lambda(x, b^i))^{\delta-1}
\right)
^{X( \Lambda(x, b^i))}
\right| X( \Lambda(x, b^i)) \ge c_4 |x|^{db^i}
\right]. 
\end{eqnarray*}
Using $1-x \le e^{-x}$, we obtain  
\begin{equation}\label{Equation: limit}
\p^k\left[
W_{z_i} \le \left(c_4 |x|^{db^i}\right)^{(1-\delta)/\beta}
\right]
~\le~
\E\left[\left.
e^{-X( \Lambda(x, b^i))^\delta}
\right| X( \Lambda(x, b^i)) \ge c_4 |x|^{db^i}
\right]
~\le~
e^{-c_4^\delta  |x|^{d\delta b^i}}.
\end{equation}
Since $|z_i-z_{i+1}| < d|x|^{b^i}$ for each $i=0, \ldots, k-1$, 
we therefore obtain
\begin{eqnarray*}
\p^k\left[ \bigcup_{i=0}^{k-1} \{ z_i \not \Leftrightarrow z_{i+1}\} \right]
&\le&
\sum_{i=0}^{k-1}
\E^k\left[
e^{-\lambda d^{-\alpha}W_{z_i}W_{z_{i+1}}|x|^{-\alpha b^i}}
\right]
\\&\le&
\sum_{i=0}^{k-1}
\E^k\left[
e^{-\lambda d^{-\alpha}W_{z_i}W_{z_{i+1}}|x|^{-\alpha b^i}}
1_{\left\{
W_{z_j} > \left(c_4 |x|^{db^j}\right)^{(1-\delta)/\beta}; \; j=i,i+1
\right\}}
\right]
+
2e^{-c_4^\delta  |x|^{d\delta b^{i+1}}}
\\&\le&
\sum_{i=0}^{k-1}
\exp\left\{
-\lambda d^{-\alpha} c_4^{2(1-\delta)/\beta} 
|x|^{b^i d(1-\delta)(1+b)/\beta}
  |x|^{-\alpha b^i}
\right\}
+
2e^{-c_4^\delta |x|^{d\delta b^{i+1}}}
\\&=&
\sum_{i=0}^{k-1}
\exp\left\{
-\lambda d^{-\alpha} c_4^{2(1-\delta)/\beta} 
|x|^{b^{i}\left(d(1-\delta)(1+b)/\beta-\alpha
\right)}
\right\}
+
2e^{-c_4^\delta |x|^{d\delta b^{i+1}}}.
\end{eqnarray*}
Since $b$ was chosen such that 
$d(1+b)/\beta > \alpha$, we can choose $\delta =\delta(b) \in (0,1)$ 
so small that $d(1-\delta)(1+b)/\beta-\alpha >0$. 
We therefore  proceed as in \eqref{Equation: similar} 
to see that there exists $m_1\ge0$ (used for the choice of $m>m_1$) so large that 
for any sufficiently large $|x|$, 
$\p^k\left[ \bigcup_{i=0}^{k-1} \{ z_i \not \Leftrightarrow z_{i+1}\} \right] 
\le \varepsilon'$. 
Using symmetry we therefore obtain for sufficiently large $|x|$,  
\begin{equation}\label{Equation: E1}
\p^k\left[  z_i \Leftrightarrow z_{i+1}
\text{ and }
z'_i \Leftrightarrow z'_{i+1} \text{ for all } i=0, \ldots ,k-1 \right]
\ge 1-2\varepsilon'.
\end{equation}
Recall that $z_0$ is the Poisson particle 
with maximal weight $W_{z_0}$ in box $\Lambda(x,1)$ with 
side length $|x|/2$ centered 
at the midpoint of the segment 
with end points $0$ and $x$. 
For every $w \ge 1$
it holds that for sufficiently large $|x|$, see also \eqref{Equation: limit} with $i=0$,  
\begin{equation*}
\p^k\left[ W_{z_0} \le w \right]
\le
\varepsilon'.
\end{equation*}
Moreover, using that box $\Lambda(x,1)$ does not 
intersect boxes $\Lambda(x, b^i)$
and  $ \Lambda'(x, b^i)$
for all $i=1, \ldots , k$, we obtain for $|x|$ 
sufficiently large, 
\begin{eqnarray*}
\p^k\left[
|\mathcal{C}(z_0)| < \infty 
\Big| W_{z_0} \ge w
\right]
&\le& 
\frac{\p\left[
|\mathcal{C}(z_0)| < \infty \Big | W_{z_0} \ge w,
\, X( \Lambda(x,1))\ge c_4|x|^d
\right]}
{\p\left[
X( \Lambda(x, b^i)) \ge c_4|x|^{db^i} \text{ and }  X( \Lambda'(x, b^i)) \ge c_4|x|^{db^i} 
\text{ for all } i=1, \ldots , k\, 
\right]}
\\&\le&
\frac{\p\left[
|\mathcal{C}(z_0)| < \infty \Big | W_{z_0} \ge w,
\, X( \Lambda(x,1))\ge c_4|x|^d
\right]}
{1-\varepsilon'},
\end{eqnarray*}
where we used \eqref{Equation: kprob} for the second inequality. 
Using Lemma \ref{Lemma:distance1} with $r=|x|/4$ we see that 
the numerator is less than $\varepsilon'$ for sufficiently large $w$ 
(note that  convergence in Lemma \ref{Lemma:distance1} 
is uniform in $r$ and $c\ge1$). 
We therefore obtain for sufficiently large $|x|$, 
\begin{equation*}
\p^k\left[
|\mathcal{C}(z_0)| = \infty 
\right]
\ge
\p^k\left[
|\mathcal{C}(z_0)| = \infty 
\Big| W_{z_0} \ge w
\right]
\p^k\left[W_{z_0} \ge w \right]
\ge
\frac{1-2\varepsilon'}{1-\varepsilon'}(1-\varepsilon')
=
1-2\varepsilon'.
\end{equation*}
Together with \eqref{Equation: E1}, this implies that 
for sufficiently large $|x|$, the event  
\begin{equation*}
E=
\{z_i \Leftrightarrow z_{i+1}
\text{ and }
z'_i \Leftrightarrow z'_{i+1} \text{ for all } i=0, \ldots ,k-1\}
\cap \{|\mathcal{C}(z_0)| = \infty\}
\end{equation*}
satisfies $\p^k[E]>1-4\varepsilon'$. 
It follows that, using \eqref{Equation: kprob}, 
\begin{eqnarray*}
\p^k[E \mid 0,x \in \mathcal{C}_\infty] 
&\ge& 
1-4\varepsilon'
/\p_{0,x}^k\left[0,x\in\mathcal{C}_\infty\right]
\\&=&
1-\frac{4\varepsilon'\p[N_k]}{
\p\left[N_k \big| 0,x\in\mathcal{C}_\infty\right]\p_{0,x}\left[0,x\in\mathcal{C}_\infty\right]}
\\&\ge&
1-\frac{4\varepsilon'}{
(1-\varepsilon'/\p_{0,x}\left[0,x\in\mathcal{C}_\infty\right])\p_{0,x}\left[0,x\in\mathcal{C}_\infty\right]}
~=~
1-\frac{4\varepsilon'}{
\p_{0,x}\left[0,x\in\mathcal{C}_\infty\right] - \varepsilon'}
~=~
1-\varepsilon'',
\end{eqnarray*}
where the last equality defines $\varepsilon''>0$. 
Note that on  event $E$, because of the choice 
of $k$ and since $z_0=z_0'$,  
\begin{eqnarray*}
d(0,x) 
&\le& 
d(0,z_k)+\sum_{i=1}^{k}\left(d(z_i,z_{i-1})+d(z'_{i-1},z'_i)\right)+d(z_k',x)
\\&=&  
d(0,z_k)+2k+d(z_k',x) 
~\le~  
d(0,z_k)+2\frac{\log \log |x|}{| \log b |}+d(z_k',x).
\end{eqnarray*}
Moreover, on  event $E$, because $z_0=z_0'$, we have that $z_k$ and 
$z'_k$ are both in the infinite connected component $\mathcal{C}_\infty$. 
If, in addition, we assume that $0 \in \mathcal{C}_\infty$, then $0$ and $z_k$ are in 
the same component $\mathcal{C}_\infty$. 
Since  $|z_k| \le 3|x|^{b^k}/4 \le m^{1/b}$
and the infinite connected component $\mathcal{C}_\infty$ 
is unique, a.s.,  it follows that on $E \cap \{0 \in \mathcal{C}_\infty\}$, 
$0$ and $z_k$ are connected within box $\Lambda_{\widetilde m}$ with 
probability arbitrarily close to $1$, for some $\widetilde m=\widetilde m(m) < \infty$. 
Note that $\widetilde m$ is independent of $x$. 
This implies for any $\kappa>0$ and sufficiently large $|x|$,
\begin{equation*}
\p^k\left[
d(0,z_k) \le \kappa \log \log |x|
\,\Big| 0 \in \mathcal{C}_\infty, \, E
\right]
\ge 
1-\varepsilon'.
\end{equation*}
By symmetry we therefore obtain for $|x|$ sufficiently large, 
\begin{equation*}
\p^k\left[
d(0,z_k)+d(z_k',x) \le 2\kappa \log \log |x|
\,\Big| 0,x \in \mathcal{C}_\infty, \, E
\right]
\ge 
1-2\varepsilon'.
\end{equation*}
Therefore, if we choose 
$\kappa = \varepsilon/(2|\log b|)$ 
and $|x|$ sufficiently large,
\begin{eqnarray*}
&&\hspace{-.5cm}
\p^k\left[\left.
d(0,x)
\le
\frac{2(1+ \varepsilon/2)\log \log |x|}{|\log b|}
\right| 0,x \in \mathcal{C}_\infty
\right]
\\&&\ge~
\p^k\left[\left.
d(0,x)
\le
\frac{2(1+ \varepsilon/2)\log \log |x|}{|\log b|}
\right| 0,x \in \mathcal{C}_\infty,\, E
\right]
\p^k\left[\left. E \right| 0,x \in \mathcal{C}_\infty\right]
\\&&\ge~
\p^k\left[\left.
d(0,z_k)+d(z_k',x)
\le
\frac{\varepsilon \log \log |x|}{|\log b|}
\right| 0,x \in \mathcal{C}_\infty,\, E
\right]
(1-\varepsilon'')
~\ge~
(1-2\varepsilon')(1-\varepsilon'').
\end{eqnarray*}
It follows that for sufficiently large $|x|$, 
using \eqref{Equation: kprob} in the last step, 
\begin{eqnarray*}
\p\left[\left.
d(0,x)
\le
\frac{2(1+ \varepsilon/2)\log \log |x|}{|\log b|}
\right| 0,x \in \mathcal{C}_\infty
\right]
&\ge&
(1-2\varepsilon')(1-\varepsilon'')
\p\left[ N_k \big | 0,x \in \mathcal{C}_\infty
\right]
\\&\ge&
(1-2\varepsilon')(1-\varepsilon'')(1-\varepsilon'/\p_{0,x}\left[0,x\in\mathcal{C}_\infty\right]).
\end{eqnarray*}
This finishes the proof of Proposition \ref{Proposition: (a)1}.
\end{Proof}

\subsection{Infinite variance of degree distribution, lower bound}

Next, we give the proof of the lower bound of statement \textit{(a)} of
Theorem \ref{Theorem: graph distance} which we recall 
in the following  proposition. 
Note that this proposition differs from the corresponding 
Theorem 5.3 of \cite{Remcoscale} in the discrete space model.

\begin{prop}\label{Proposition: (a)2}
Assume  $\alpha>  d$ 
and $\tau=\beta \alpha/d \in (1, 2)$. 
For every $\lambda > \lambda_c=0$ there exists $\eta_1 > 0$ 
such that 
\begin{equation*}
\lim_{|x| \rightarrow \infty}
\p_{0,x}\left[
d(0,x)
\ge 
\eta_1
\frac{2\log \log |x|}{|\log\kappa|}
\right]
=1,
\end{equation*}
with $\kappa= \alpha(\beta \wedge 1)/d -1 \in (0,1)$.
\end{prop}

~

\begin{Proof}[of Proposition \ref{Proposition: (a)2}]
We modify the proof of Theorem 5.3 of \cite{Remcoscale} to our situation. 
Choose $\vartheta>1$ and $\mu>0$ such that 
\begin{equation*}
d/\vartheta-d\kappa+\mu<0 
\quad \text{and}\quad
\mu<d\kappa.
\end{equation*}
Note that this choice is possible  
since the above constraints require $\vartheta>1/\kappa$ and 
$\mu\in(0,d(\kappa-1/\vartheta))$. 
For $x\in X$ and $n\in \N$ we define the random variable 
\begin{equation*}
S_n(x) = \sup_{y\in X:\, d(x,y)\le n} |x-y|,
\end{equation*}
which represents the Euclidean distance between  $x$ 
and the furthest particle that can be reached from 
$x$ using at most $n$ edges. 
For $r>0$ we denote by $ B(r)$ the ball of 
(Euclidean) radius $r$ around the origin. 
For $t>1$ we obtain, using  $1-e^{-x} \le x\wedge 1$, 
\begin{eqnarray*}
\p_{0,x}\left[ S_{n-1}(0) < t^{1/\vartheta}, \;S_n(0) \ge t \right]
&\le&
\p_{0,x}\left[ 
\exists z \in  B(t^{1/\vartheta})\cap X, \,  z' \in  B(t)^c \cap X
\text{ such that } z \Leftrightarrow z'
\right]
\\&\le&
\E_{0,x}\left[
\sum_{z \in  B(t^{1/\vartheta})\cap X}
\sum_{  z' \in  B(t)^c\cap X}
\E\left[\left.  \frac{\lambda W_zW_{z'}}{|z-z'|^\alpha} \wedge 1 \right | X \right]
\right].
\end{eqnarray*}
For two i.i.d.~random variables $W_1$ and $W_2$ 
having a Pareto distribution with scale parameter 
$1$ and shape parameter $\beta$ we obtain 
for $u \ge 1$, using integration by parts in the first step, 
\begin{eqnarray*}
&&\hspace{-.75cm}
\E\left[ \frac{ W_1W_2}{u} \wedge 1 \right]
~=~
\frac{1}{u}
+
\frac{1}{u}  \int_1^u \p[W_1W_2> v]dv
~=~
\frac{1}{u} 
+ \frac{1}{u} \int_1^u v^{-\beta}(1+ \beta \log v) dv
\\&&\le~
(1+ \beta \log u)\left( u^{-(\beta \wedge 1)} 
+ \frac{1}{u} \int_1^u v^{-\beta}dv\right)
~\le~
\max\{1+\log u,1+1_{\{\beta\neq1\}}/|\beta-1|\} \left( 1+ \beta \log u\right ) u^{-(\beta \wedge 1)},
\end{eqnarray*}
where the last step follows by distinguishing between the cases 
$\beta=1$, $\beta > 1$ and $\beta < 1$. 
This provides  for $u\ge 1$,
\begin{equation}\label{Equation: Pareto}
\E\left[ \frac{ W_1W_2}{u} \wedge 1 \right]
\le 
(1+1_{\{\beta\neq1\}}/|\beta-1|) \left( 1+ \max\{1,\beta \} \log u\right )^2 u^{-(\beta \wedge 1)}.
\end{equation}
Choose $t$ so large that $\lambda^{-1}(t-t^{1/\vartheta})^\alpha\ge 1$ 
which together with \eqref{Equation: Pareto} implies that 
\begin{eqnarray*}
&&\hspace{-.5cm}
\p_{0,x}\left[ S_{n-1}(0) < t^{1/\vartheta}, \;S_n(0) \ge t \right]
~\le~
\E_{0,x}\left[
\sum_{z \in  B(t^{1/\vartheta})\cap X}
\sum_{  z' \in  B(t)^c\cap X}
\E\left[ \left. \frac{ W_zW_{z'}}{\lambda^{-1}|z-z'|^\alpha} \wedge 1 \right | X\right]
\right]
\\&&\hspace{-.3cm}\le
(1+1_{\{\beta\neq1\}}/|\beta-1|) \E_{0,x}\left[
\sum_{z \in  B(t^{1/\vartheta})\cap X}
\sum_{  z' \in  B(t)^c\cap X}
\left( 1+ \max\{1,\beta \} \log \left( \lambda^{-1} |z-z'|^{\alpha} \right)\right )^2 
\left(\lambda^{-1} |z-z'|^{\alpha}\right)^{-(\beta \wedge 1)}
\right].
\end{eqnarray*}
Choose $t$ so large that  
$(1+1_{\{\beta\neq1\}}/|\beta-1|)
\left( 1+ \max\{1,\beta \} \log \left( \lambda^{-1} |z-z'|^{\alpha} \right)\right )^2
\lambda^{\beta \wedge 1}
\le 
|z-z'|^{\mu}$ 
for all $z \in  B(t^{1/\vartheta})$ and $z' \in  B(t)^c$. 
It follows that for sufficiently large $t$, 
note that $d(\kappa+1)= \alpha(\beta \wedge 1)$, 
\begin{eqnarray*}
\p_{0,x}\left[ S_{n-1}(0) < t^{1/\vartheta}, \;S_n(0) \ge t \right]
&\le&
\E_{0,x}\left[ 
\sum_{z \in  B(t^{1/\vartheta})\cap X}
\sum_{  z' \in  B(t)^c\cap X}
|z-z'|^{-\alpha(\beta \wedge 1) + \mu}
\right]
\\&=&
\E_{0,x}\left[ 
\sum_{z \in  B(t^{1/\vartheta})\cap X}
\E_{0,x}\left[
\sum_{  z' \in  B(t)^c\cap X}
|z-z'|^{-d(\kappa+1) + \mu}
\right]
\right].
\end{eqnarray*}
We estimate the right-hand side under the unconditional 
measure $\p$ instead of $\p_{0,x}$. Note that the 
tail behavior is the same under both measures. 
We obtain  
\begin{eqnarray*}
\E\left[ 
\sum_{z \in  B(t^{1/\vartheta})\cap X}
\E\left[
\sum_{  z' \in  B(t)^c\cap X}
|z-z'|^{-d(\kappa+1) + \mu}
\right]
\right]
&=&
\E\left[ 
\sum_{z \in  B(t^{1/\vartheta})\cap X}
\nu \int_{z' \in  B(t)^c}
|z-z'|^{-d(\kappa+1) + \mu}dz'
\right]
\\&\le&
\E\left[ 
\sum_{z \in  B(t^{1/\vartheta})\cap X}
\nu \int_{|z'| \ge t-t^{1/\vartheta}}
|z'|^{-d(\kappa+1) + \mu}dz'
\right]
\\&=&
\nu^2 v_d t^{d/\vartheta}
\int_{|z'| \ge t-t^{1/\vartheta}}
|z'|^{-d(\kappa+1) + \mu}dz'.
\end{eqnarray*}
We therefore obtain for an appropriate constant $c_5>0$ and 
for $t$ sufficiently large, 
\begin{equation}\label{Equation: boundA}
\p_{0,x}\left[ S_{n-1}(0) < t^{1/\vartheta}, \;S_n(0) \ge t \right]
\le
c_5 t^{d/\vartheta-d\kappa + \mu}.
\end{equation}
%
%
Define $f:\N_0^2 \rightarrow (0,\infty)$ by $f(m,n)=m^{\vartheta^n}$ for all $m,n\in\N_0$. 
Observe 
\begin{itemize}[leftmargin=15pt]
\item[1)] 
$\sum_{k=2}^{\infty}f(2,k)^{d/\vartheta-d\kappa + \mu}<\infty$  
	because $\vartheta>1$ and $d/\vartheta-d\kappa + \mu<0$;
\item[2)] 
for all $m\ge2$ and sufficiently small $\eta_1=\eta_1(m)>0$: 
         $f\left(m,\left \lceil \eta_1\frac{\log \log |x|}{|\log\kappa|} \right \rceil \right)\le|x|/2$ 
         for all sufficiently large $|x|$. 
\end{itemize}
We choose $m_0\ge2$ so large that \eqref{Equation: boundA}  
holds true for all $t=f(m,n)$ with $m\ge m_0$ and $n\ge2$. 
Using \eqref{Equation: boundA} and induction we obtain  
for each $n\ge2$ and $m\ge m_0$, 
note that $f(m,n)^{1/\vartheta}=f(m,n-1)$, 
\begin{eqnarray*}
\p_{0,x}\left[ S_n(0) \ge f(m,n) \right]
&\le&
\p_{0,x}\left[ S_{n-1}(0) \ge f(m,n-1) \right]
+
\p_{0,x}\left[ S_{n-1}(0) < f(m,n-1), \;S_n(0) \ge f(m,n) \right]
\\&\le&
\p_{0,x}\left[ S_{n-1}(0) \ge f(m,n-1) \right]
+
c_5 
f(m,n)^{d/\vartheta-d\kappa + \mu}
\\&\le&
\p_{0,x}\left[ S_{1}(0) \ge f(m,1) \right]
+
c_5 
\sum_{k=2}^{n}f(m,k)^{d/\vartheta-d\kappa + \mu}
\\&\le&
\p_{0,x}\left[ \exists\; y \in X \text{ with } |y| \ge  f(m,1) \text{ and } 0 \Leftrightarrow y \right]
+
c_5 
 \sum_{k=2}^{\infty}f(m,k)^{d/\vartheta-d\kappa + \mu}.
\end{eqnarray*}
Note that the right-hand side is independent of $n\ge2$ and is 
finite for any $m\ge2$. 
Since $f(m,k)^{d/\vartheta-d\kappa + \mu}$ is decreasing 
in $m\ge2$, there exists $m\ge m_0$ such that 
the right-hand side is less than $\varepsilon$ for fixed $\varepsilon>0$.  
We 
finally obtain for sufficiently small $\eta_1=\eta_1(m)>0$ 
and for all sufficiently large $|x|$, 
set $n(x)= \left\lceil \eta_1\frac{\log \log |x|}{|\log\kappa|} \right \rceil \ge 2$, 
\begin{eqnarray*}
\p_{0,x}\left[
d(0,x)
\le
\eta_1
\frac{2\log \log |x|}{|\log\kappa|}
\right]
&\le&
\p_{0,x}\left[d(0,x)\le 2n(x)\right]
\\&\le& 
\p_{0,x}\left[S_{n(x)}(x)\ge |x|/2\right]
+
\p_{0,x}\left[S_{n(x)}(0)\ge |x|/2\right]
\\&=&
2
\p_{0,x}\left[
S_{n(x)}(0)
\ge 
|x|/2
\right]
\\&\le&
2
\p_{0,x}\left[
S_{n(x)}(0)
\ge 
f\left(m,\left\lceil \eta_1\frac{\log \log |x|}{|\log\kappa|} \right \rceil \right)
\right]
~\le~
2\varepsilon,
\end{eqnarray*}
which finishes the proof of Proposition \ref{Proposition: (a)2}.
\end{Proof}

\subsection{Finite variance of degree distribution, case 1, lower bound}

In the following we give the proof of 
part \textit{(b1)} of Theorem \ref{Theorem: graph distance}. 
We first prove the lower bound which 
follows from the following proposition.

\begin{prop}\label{Proposition: (b1)1}
Assume $\alpha> d$ and 
$\tau=\beta \alpha/d > 2$. For every $\lambda > \lambda_c$ 
there exists $\eta' > 0$ such that 
\begin{equation*}
\p_{0,x}\left[
d(0,x) \ge \eta' \log|x|
\right]=1.
\end{equation*}
\end{prop}

~

\begin{Proof}[of Proposition \ref{Proposition: (b1)1}]
Choose $n \in \N$, $0,x\in X$ and set $x_0=0$ and $x_n=x$. 
As in \eqref{Equation: g} we obtain, 
the first sum is over all self-avoiding paths 
of length $n$ starting from $0$, note that $x_n=x$ is now fixed, 
\begin{equation*}
\p_{0,x}\left[
d(0,x)=n
\right]
\le
\E_{0,x}\left[
\sum_{\text{$(x_1, \ldots, x_{n})$ s.a.}}
\prod_{i=1}^n p_{x_{i-1}x_i}
\right]
\le
\nu^{n-1}
\int_{\R^d}\cdots\int_{\R^d}
\prod_{i=1}^{n}
 h(x_{i}-x_{i-1})
dx_1 \cdots dx_{n-1}, 
\end{equation*}
where for $y\in\R^d$ we define function $h$ by, 
recall function $g$ defined in \eqref{Equation: function g},
\begin{equation*}
 h(y) 
~=~ 
 g(\lambda^{-1}|y|^\alpha)
 ~=~
1_{\left\{ |y| < \lambda^{1/\alpha} \right\}} 
+
1_{\left\{ |y| \ge \lambda^{1/\alpha} \right\}} 
c_2\lambda^{(\beta/2 \wedge 1)} \left( 1+ \max\{2,\beta \}\log(\lambda^{-1}|y|^\alpha)\right) |y|^{-\alpha(\beta/2 \wedge 1)}.
\end{equation*}
Note that $h$ is integrable because $\alpha > d$ and $\tau=\beta \alpha/d > 2$. 
Using $x_0=0$ and $x_n=x$, and 
substituting inductively $x_i$ by $x_i-\sum_{l=1}^{i-1}x_l$ for $i=1,\ldots, n-1$, it follows that 
\begin{equation*}
\p_{0,x}\left[
d(0,x)=n
\right]
~\le~
\nu^{n-1}
\int_{\R^d}\cdots\int_{\R^d} 
\left(\prod_{i=1}^{n-1}h(x_i)\right)h\left(x-\sum_{i=1}^{n-1}x_i\right)dx_1 \cdots dx_{n-1}.
\end{equation*}
We condition on $1_{\{|x_i| < |x|/n\}}$ and $1_{\{|x_i| \ge |x|/n\}}$ for all 
$i=1,\ldots,n-1$. 
Note that if $|x_i| < |x|/n$ for all $i=1, \ldots , n-1$ we have  
$|x- \sum_{i=1}^{n-1}x_i| \ge |x|/n$, and we bound the corresponding 
factor in the integral by $\sup_{ y \in \R^d: \, |y| \ge |x|/n} h(y)$. 
Otherwise, by exchangeability of the $x_i$'s, there are $n-1$ different cases where 
at least one of the $x_i$'s satisfies $|x_i| \ge |x|/n$. 
In each of these $n-1$ cases we bound one corresponding factor in the integral 
by $\sup_{ y \in \R^d: \, |y| \ge |x|/n} h(y)$. Note that the  
restriction on $x- \sum_{i=1}^{n-1}x_i$ then drops and we obtain 
\begin{equation}\label{Equation: Convolution}
\p_{0,x}\left[
d(0,x)=n
\right]
\le 
 n \left( \sup_{ y \in \R^d: \, |y| \ge |x|/n} h(y) \right) 
\left(\nu \int_{\R^d}h(y) dy \right )^{n-1},
\end{equation} 
where $ \nu\int_{\R^d}h(y) dy < \infty$ since $h$ is integrable. 
Next, we bound the supremum on the right-hand side of \eqref{Equation: Convolution}. 
Choose $\eta >0$ and let $|x|$ be so large that $\eta \log|x| \ge 1$. 
Choose $n \in \N$ with $n \le \eta \log|x|$. 
Let $\mu \in (0, \alpha(\beta/2 \wedge 1))$ and choose 
$|x|$ so large that any $y\in \R^d$ with $|y| \ge |x|/n$ 
satisfies 
\begin{equation*}
c_2\lambda^{(\beta/2 \wedge 1)}
\left(1+ \max\{2,\beta \}\log(\lambda^{-1}|y|^\alpha)\right) \le |y|^\mu.
\end{equation*}
If, in addition, $|x|$ is so large that 
$|x| / n > \lambda^{1/\alpha}$, then for any $y\in \R^d$ with $|y| \ge |x|/n> \lambda^{1/\alpha}$,
\begin{equation*}
h(y) 
= 
c_2\lambda^{(\beta/2 \wedge 1)} \left( 1+ \max\{2,\beta \}\log(\lambda^{-1}|y|^\alpha)\right) |y|^{-\alpha(\beta/2 \wedge 1)}
\le
|y|^{-\alpha (\beta/2 \wedge 1)+\mu}
\le
n^{\alpha (\beta/2 \wedge 1)-\mu}
|x|^{-\alpha (\beta/2 \wedge 1)+\mu}.
\end{equation*} 
We finally obtain for all $\eta >0$ and  $1\le n \le \eta \log|x|$ with $|x|$ sufficiently large,
\begin{eqnarray*}
\sup_{ y \in \R^d: \, |y| \ge |x|/n} h(y)
\le 
\eta^{\alpha (\beta/2 \wedge 1)-\mu} (\log|x|)^{\alpha (\beta/2 \wedge 1)-\mu}
|x|^{-\alpha (\beta/2 \wedge 1)+\mu}.
\end{eqnarray*} 
Together with \eqref{Equation: Convolution} we obtain for any $\eta>0$ and  
$1\le n \le \eta \log|x|$ with $|x|$ sufficiently large, 
\begin{eqnarray*}
\p_{0,x}\left[
d(0,x)=n
\right]
&\le&
 n
\left(\nu \int_{\R^d}h(y) dy \right )^{n-1}
\eta^{\alpha (\beta/2 \wedge 1)-\mu} (\log|x|)^{\alpha (\beta/2 \wedge 1)-\mu}
|x|^{-\alpha (\beta/2 \wedge 1)+\mu}
\\&\le& 
 \left( 1+\nu\int_{\R^d}h(y) dy \right )^{\eta \log|x|}
\eta^{\alpha (\beta/2 \wedge 1)-\mu+1} (\log|x|)^{\alpha (\beta/2 \wedge 1)-\mu+1}
|x|^{-\alpha (\beta/2 \wedge 1)+\mu}
\\&=& 
\eta^{\alpha (\beta/2 \wedge 1)-\mu+1} (\log|x|)^{\alpha (\beta/2 \wedge 1)-\mu+1}
|x|^{-\alpha (\beta/2 \wedge 1)+\mu+\eta \log \left(1+ \nu\int_{\R^d}h(y) dy \right )}
~\le~
|x|^{-\delta},
\end{eqnarray*} 
where the last inequality holds for some $\delta>0$ whenever 
$|x|$ is sufficiently large and $\eta>0$ is chosen so small that 
$-\alpha (\beta/2 \wedge 1)+\mu+\eta \log \left(1+ \nu\int_{\R^d}h(y) dy \right )<0$. 
We conclude that there exist $\eta'>0$ and $\delta>0$ 
such that for all $|x|$ sufficiently large, 
\begin{equation*}
\p_{0,x}\left[
d(0,x) \le \eta' \log|x|
\right]
=
\sum_{1\le n \le \eta' \log |x|} \p_{0,x}\left[d(0,x)=n\right]
\le 
\eta' (\log|x|)
|x|^{-\delta},
\end{equation*}
which converges to $0$ as $|x| \to \infty$. 
\end{Proof}

\subsection{Finite variance of degree distribution, case 1, upper bound}

In order to finish the proof of 
statement \textit{(b1)} of Theorem \ref{Theorem: graph distance} 
it remains to show the corresponding upper bound on the graph distances. 
The result  follows from the following proposition, 
see also Proposition 4.1 of \cite{biskup}.

\begin{prop}\label{Proposition: 4.1}
Assume $\alpha \in (d,2d)$ and $\tau=\beta \alpha/d > 2$, 
and choose $\lambda > \lambda_c$. 
For each $\varepsilon>0$ and $\Delta' > \Delta= \log(2)/\log(2d/\alpha)$ 
there exists $N_0 < \infty$ such that for all 
$x,y \in \R^d$ with $|x-y| \ge N_0$, 
\begin{equation*}
\p_{x,y}\left[
d(x,y) \ge (\log|x-y|)^{\Delta'}
,\, x,y \in \mathcal{C}_\infty
\right]\le \varepsilon.
\end{equation*}
\end{prop}

Note that the latter statement implies 
\begin{eqnarray*}
\p\left[ 
d(x,y) \le (\log|x-y|)^{\Delta'}
\Big | x,y \in \mathcal{C}_\infty
\right]
&=& 
1-\frac{\p_{x,y}\left[ 
d(x,y) > (\log|x-y|)^{\Delta'}
,\,  x,y \in \mathcal{C}_\infty
\right]}{\p_{x,y}\left[ x,y \in \mathcal{C}_\infty \right]}
\\&\ge& 
1-\frac{\varepsilon}{\p_{x,y}\left[ x,y \in \mathcal{C}_\infty \right]}. 
\end{eqnarray*}

To prove Proposition \ref{Proposition: 4.1} we use 
the concept of hierarchies of particles. 
For $k\ge1$ we call an element $\sigma \in \{0,1\}^k$, such 
as $\sigma = 01110001$, a hierarchical 
index. 
If $k=0$, $\sigma \in \{0,1\}^k$ denotes the empty string. 
For $\sigma_1\in \{0,1\}^k$ and $\sigma_2 \in \{0,1\}^l$, 
$k,l\ge 1$, we denote by $\sigma_1\sigma_2$ the concatenation 
of $\sigma_1$ and $\sigma_2$. 
Then \cite{biskup} provides the following definition 
of a hierarchy. 

\begin{defi}\label{Definition: hierarchy}
For $m \ge 1$ and two distinct particles $x,y \in X$ 
we say that the set of particles 
\begin{equation*}
\mathcal{H}_m(x,y) 
= 
\left\{
z_\sigma \in X \, \Big| \sigma \in \{0,1\}^k , \, k=1, \ldots , m
\right\} \subset X
\end{equation*}
is a hierarchy of depth $m$ connecting $x$ and $y$ if 
\begin{itemize}
\item[1.]
$z_0=x$ and $z_1=y$;
\item[2.]
$z_{\sigma 00}= z_{\sigma 0}$ and $z_{\sigma 11}= z_{\sigma 1}$ 
for all $\sigma \in \{0,1\}^k$ and $k=0, \ldots, m-2$;
\item[3.]
for all $\sigma \in \{0,1\}^k$ and $k=0, \ldots,m-2$ there is an edge 
between $z_{\sigma01}$ and $z_{\sigma10}$ as 
long as $z_{\sigma01} \neq z_{\sigma10}$;
\item[4.]
each edge as in 3.~appears only once in $\mathcal{H}_m(x,y)$.
\end{itemize}
For $\sigma \in \{0,1\}^{m-2}$ 
we call the pairs of particles $(z_{\sigma00},z_{\sigma01})$ and 
$(z_{\sigma11},z_{\sigma10})$  ``gaps''.
\end{defi}

\begin{figure}
\begin{center}
\includegraphics[width=\textwidth]{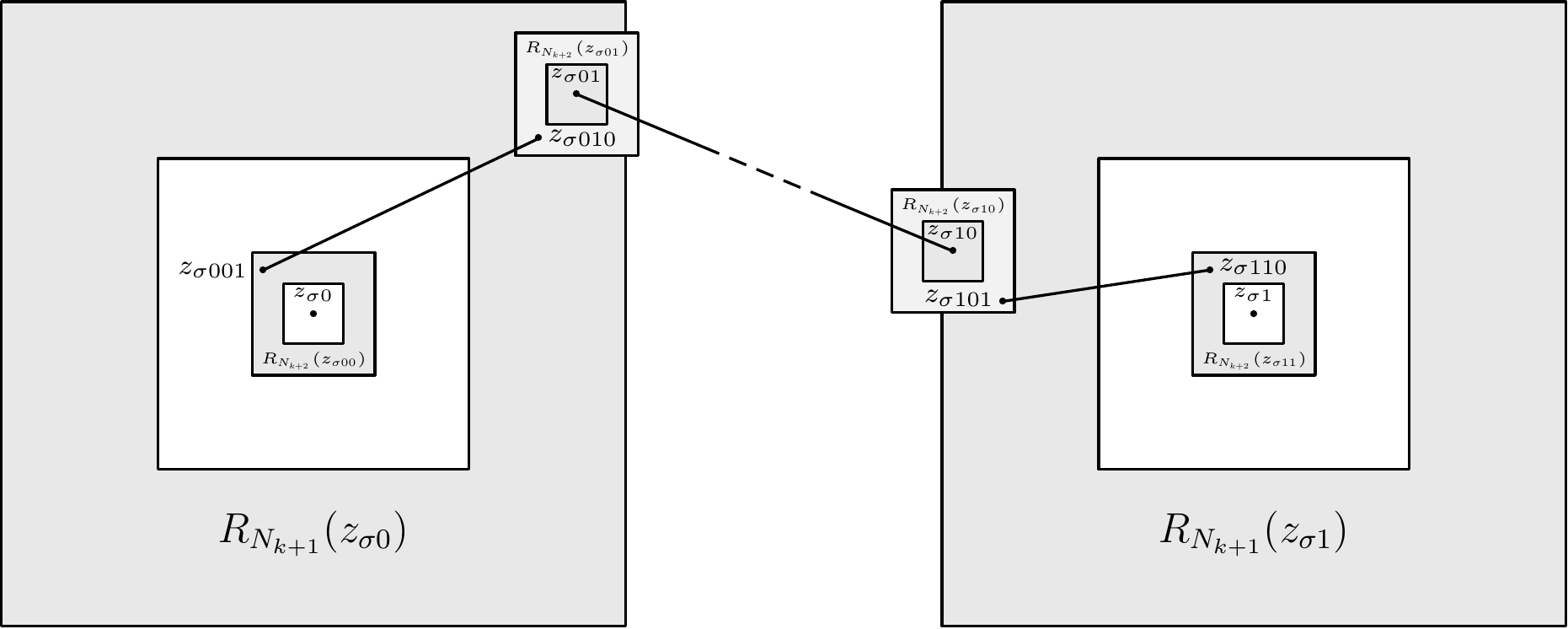}
 \end{center}
\caption{\footnotesize
Illustration of a hierarchy in $\mathcal{B}_m$. 
Assume we are given $z_{\sigma0}$ and $z_{\sigma1}$ 
for some $\sigma \in \{0,1\}^k$ and $k=0, \ldots , m-3$. 
We consider an edge between $R_{N_{k+1}}(z_{\sigma 0})$ and 
$R_{N_{k+1}}(z_{\sigma 1})$ and call the end points 
$z_{\sigma 01}$ and $z_{\sigma 10}$, respectively.  
Then, set $z_{\sigma00}=z_{\sigma0}$ and consider an edge 
between $R_{N_{k+2}}(z_{\sigma 00})$ and 
$R_{N_{k+2}}(z_{\sigma 01})$ and call the end points  $z_{\sigma 001}$ and 
$z_{\sigma 010}$, respectively. 
Similarly, we define $z_{\sigma 110}$ and 
$z_{\sigma 101}$. 
To get a path of edges from $z_{\sigma0}$ to $z_{\sigma1}$ 
it remains to connect the two particles of each gap $(z,z')$. 
(If we assume  $m=3$, the pairs 
$(z_{\sigma0},z_{\sigma001})$, 
$(z_{\sigma01},z_{\sigma010})$, 
$(z_{\sigma11},z_{\sigma101})$ and 
$(z_{\sigma1},z_{\sigma110})$ 
in the figure are the gaps of the illustrated hierarchy in $\mathcal{B}_3$.)
Such particles are likely connected by a short path of edges
because they are $(\rho,\ell)$-dense and close to each other if $m$ is sufficiently large,  
conditional on $\mathcal{B}_m$. 
}
\label{Figure1} 
\end{figure}

Recall that for $x \in \R^d$ and $n\in(0,\infty)$ we write 
$\Lambda_n(x)=x +[-n,n)^d$ for the box with center $x$ and 
side length $2n$,  
and for $x \in X$ we write $\mathcal{C}_n(x)$ for the set of 
particles in $\Lambda_n(x) \cap X$ that are connected 
to $x$ within $\Lambda_n(x)$. For any $L>0$ and $x\in\R^d$ we set 
annulus
$R_L(x)= \Lambda_{L}(x) \setminus \Lambda_{L/2}(x)$.
Moreover, for $\ell\in(0,L)$ and $\rho \in (0,1)$ we define by 
\begin{equation*}
\mathcal{D}_L^{(\rho,\ell)}(x) 
= 
\left\{
z \in R_L(x)\cap X \Big | |\mathcal{C}_\ell(z)| \ge \rho (2\ell)^d 
\right\},
\end{equation*}
the set of $(\rho,\ell)$-dense particles in $R_L(x)$. 
Note that this 
definition differs from $\mathcal{D}_L^{(\rho,\ell)}$ defined in Section \ref{Results:Perc} because 
we now consider the particles in annulus $R_L(x)$ instead of the 
particles in box $\Lambda_L$. 
For $x,y\in\R^d$, $\alpha\in(d,2d)$ and $\gamma \in (\alpha/(2d),1)$ 
we define for $m\ge1$,
\begin{equation*}
N_m=N^{\gamma^m} \quad \text{with $N=|x-y|$}.
\end{equation*}
We denote by $\mathcal{B}_m=\mathcal{B}_{m,\gamma}^{(\rho,\ell)}(x,y)$ 
the event that there exists a hierarchy $\mathcal{H}_m(x,y)$ of depth $m$ 
connecting $x$ and $y$ such that  
for all $k=0, \ldots , m-2$ and all $\sigma \in \{0,1\}^k$, 
\begin{equation*}
z_{\sigma01} \in \mathcal{D}_{N_{k+1}}^{(\rho,\ell)}(z_{\sigma0})
\quad \text{and} \quad 
z_{\sigma10} \in \mathcal{D}_{N_{k+1}}^{(\rho,\ell)}(z_{\sigma1}),
\end{equation*}
see also (4.5) of \cite{biskup} and Figure \ref{Figure1}  for an illustration. 
Moreover, we denote by $\mathcal{T}=\mathcal{T}^{(\rho,\ell)}(x,y)$ 
the event that $x$ and $y$ are $(\rho,\ell)$-dense. 
Note that the event $\mathcal{B}_m \cap \mathcal{T}$ 
ensures that there is a hierarchy $\mathcal{H}_m(x,y)$ of depth $m$ 
connecting $x$ and $y$ such that
all particles in this hierarchy are $(\rho,\ell)$-dense and lie 
in the corresponding annulus $R_{N_k}(z_\sigma)$. 
Finally, given $\mathcal{B}_m$, we denote by $\mathcal{S}=\mathcal{S}^{(\ell)}$ 
the event that all gaps $(z,z')$ 
in a hierarchy in $\mathcal{B}_m$ satisfy 
\begin{equation}\label{Equation: S}
X( \Lambda_\ell(z) ) \le e\nu (2\ell)^d
\quad \text{and} \quad 
X( \Lambda_\ell(z') ) \le e\nu (2\ell)^d.
\end{equation}
The event $\mathcal{S}$ ensures that we do not have too 
many particles in  boxes $\Lambda_\ell(z)$ and $\Lambda_\ell(z')$,  
in particular, the graph distance of connected paths within such boxes is 
bounded by $e\nu(2\ell)^d$.
%
%
%
%
For the proof of Proposition \ref{Proposition: 4.1} we use 
Lemma \ref{Lemma: 4.2,4.3} below. 
Part \textit{(i)} of Lemma \ref{Lemma: 4.2,4.3} corresponds 
to Lemma 4.2 of \cite{biskup}. The only difference lies in the additional 
event $\mathcal{S}$ defined by \eqref{Equation: S}. 
Part \textit{(ii)} of Lemma \ref{Lemma: 4.2,4.3} is the continuum space 
analogue to Lemma 4.3 of \cite{biskup} and  
it shows that the event $\mathcal{B}_m$ occurs with  sufficiently high probability. 

\begin{lemma}\label{Lemma: 4.2,4.3}
Choose $\alpha \in (d,2d)$. 
For all $\varepsilon \in (0,1)$, $\gamma \in (\alpha/(2d),1)$,  
$\Delta' > \log(2)/\log(1/\gamma)$ and $\alpha' \in (\alpha, 2d\gamma)$  
there exist $N' = N'(\varepsilon,\gamma, \Delta') < \infty$, 
$\rho_0 \in (0,1)$ and  a constant $c_6>0$ such that the following holds true:
for all $x,y \in \R^d$ with $N=|x-y| \ge N'$ let $m\in\N$ be the maximal integer such that
\begin{equation}\label{Equation: (4.7)}
m\log(1/\gamma) \le \log \log N - \varepsilon \log \log \log N.
\end{equation}
For all $\rho \in (0,\rho_0)$ and  $\ell\in(N_m,2N_m)$ in the definitions 
of $\mathcal{B}_m$, $\mathcal{T}$ and $\mathcal{S}$ we obtain
\begin{enumerate}
\item[(i)]
$
\p_{x,y}\left[
\left\{ d(x,y) \ge (\log N)^{\Delta'}\right\}
\cap \mathcal{B}_m \cap \mathcal{T} \cap \mathcal{S}
\right] \le \varepsilon;
$
\item[(ii)]
$
\p_{x,y}\left[
\mathcal{B}_m^c
\right]
\le 
2^{m+1}e^{ -c_6N_m^{2d\gamma-\alpha'}};
$
\item[(iii)]
$
\p_{x,y}\left[ \mathcal{B}_m \cap\mathcal{S}^c \right]
\le
\varepsilon.
$
\end{enumerate}
\end{lemma}

Note that choice \eqref{Equation: (4.7)} implies, see also (4.9) of \cite{biskup},
\begin{equation}\label{Equation: (4.9)}
2^m \le (\log N)^{\log2/\log(1/\gamma)}
\quad \text{and} \quad 
e^{(\log \log N)^\varepsilon} \le N_m \le e^{(1/\gamma)(\log \log N)^\varepsilon}.
\end{equation}

~

\begin{Proof}[of Lemma \ref{Lemma: 4.2,4.3}]  
We first prove \textit{(i)}. Let $\rho\in(0,1)$. Using the additional event $\mathcal{S}$ we 
argue as in the proof of Lemma 4.2 of \cite{biskup} 
to see that it remains to estimate the probability of 
the event $\mathcal{A}_m \cap \mathcal{B}_m \cap \mathcal{T} \cap \mathcal{S}
\subset\mathcal{A}_m \cap \mathcal{B}_m \cap \mathcal{T}$ 
for $N$ sufficiently large, where $\mathcal{A}_m$ is the event that for any hierarchy 
in $\mathcal{B}_m$ there exists a gap $(z,z')$ 
such that there is no edge between
the sets $\mathcal{C}_\ell(z)$ 
and $\mathcal{C}_\ell(z')$. 
Let $z\in\R^d$ be such that the event 
$\mathcal{A}_m \cap \mathcal{B}_m \cap \mathcal{T}$ only depends 
on the particles in $\Lambda_N(z)$ and edges with end points 
in $\Lambda_N(z)$ (which exists if $N$ is sufficiently large). 
For $l\ge1$ and $x=x_0,y=x_1,x_2, \ldots , x_l \in \Lambda_N(z)$ 
we consider the edge probabilities 
\begin{equation*}
\widetilde p_{x_ix_j} ~=~1-\exp\{-\lambda|x_i-x_j|^{-\alpha}\}
~\le~1-\exp\{-\lambda W_{x_i}W_{x_j}|x_i-x_j|^{-\alpha}\}~=~p_{x_ix_j},
\end{equation*}
a.s., for every $i \neq j \in \{0,\ldots, l\}$, and denote by 
$\widetilde\p_X$ the probability measure of the resulting 
edge configurations. 
We obtain 
\begin{equation*}
\p_{x,y}\Big[ 
\mathcal{A}_m \cap \mathcal{B}_m \cap \mathcal{T}
\;\Big |\; X\cap\Lambda_N(z)=\{x,y,x_1, \ldots, x_l\}
\Big]
~\le~
\widetilde\p_X\left[ 
\mathcal{A}_m \cap \mathcal{B}_m \cap \mathcal{T}
\right].
\end{equation*}
Using the same arguments as in \cite{biskup} we obtain 
for all $N$  sufficiently large, 
\begin{equation}\label{Equation: trick}
\widetilde\p_X\left[ 
\mathcal{A}_m \cap \mathcal{B}_m \cap \mathcal{T}
\right]
~\le~
\varepsilon,
\end{equation}
note that the left-hand side is zero for $l$ too small. 
This proves $(i)$ by additionally integrating over the Poisson cloud restricted 
to $\Lambda_N(z)$ 
(note that the right-hand side of \eqref{Equation: trick} does not depend on the Poisson cloud).

~

To prove \textit{(ii)} we note that 
\begin{equation*}
\p_{x,y}\left[
\mathcal{B}_m^c
\right]
\le
\p_{x,y}\left[\mathcal{B}_{m-1} \cap \mathcal{B}_{m}^c \right] + \p_{x,y}\left[\mathcal{B}_{m-1}^c\right]
\le
\sum_{k=1}^{m-1}
\p_{x,y}\left[\mathcal{B}_{k} \cap \mathcal{B}_{k+1}^c \right] + \p_{x,y}\left[\mathcal{B}_{1}^c\right]
=
\sum_{k=1}^{m-1}
\p_{x,y}\left[\mathcal{B}_{k} \cap \mathcal{B}_{k+1}^c \right].
\end{equation*}
Hence, it is sufficient to show that  
there exist $\rho_0 \in (0,1)$ and  a constant $c_6>0$ such that 
for all $\rho \in (0,\rho_0)$ and for $\ell\in(N_m,2N_m)$ in the definition 
of $\mathcal{B}_m$ we have  
for sufficiently large $N$, 
\begin{equation*}
\p_{x,y}\left[\mathcal{B}_k \cap \mathcal{B}_{k+1}^c \right] 
\le 2^{k+1} e^{ -c_6N_m^{2d\gamma-\alpha'} }, 
\end{equation*}
for all $k=1, \ldots, m-1$, with $m$ as in \eqref{Equation: (4.7)}. For $k=1, \ldots, m-1$
we denote by $\mathcal{B}_k'$ the event that there is a hierarchy 
$\mathcal{H}_k(x,y)$ of depth $k$ connecting $x$ and $y$ such that 
for each $j=0, \ldots, k-2$ and each $\sigma \in \{0,1\}^j$,
\begin{equation*}
z_{\sigma01} \in R_{N_{j+1}}(z_{\sigma0})
\quad \text{and} \quad 
z_{\sigma10} \in R_{N_{j+1}}(z_{\sigma1}).
\end{equation*}
By definition we obtain $\mathcal{B}_k \subset \mathcal{B}_k'$. 
For all $\ell\in(N_m,2N_m)$ and 
$ \rho \in (0,\rho_0)$, where $\rho_0$ will be chosen below, we define the events
\begin{eqnarray*}
\mathcal{A}_1 &=& \big\{
\text{for any hierarchy in $\mathcal{B}_k'$ there exists $\sigma \in \{0,1\}^k$ 
such that $|\mathcal{D}_{N_k}^{(\rho,\ell)}(z_\sigma)| \le  \rho_0 (2N_k)^d$}
\big\},
\\
\mathcal{A}_2 &=& \big\{
\text{for any hierarchy in $\mathcal{B}_k'$ there exists a gap $(z,z')$ 
such that there is no edge between} \\
&&\hspace{8.5cm}
\text{the sets $\mathcal{D}_{N_k}^{(\rho,\ell)}(z)$ 
and $\mathcal{D}_{N_k}^{(\rho,\ell)}(z')$}
\big\}.
\end{eqnarray*}
$\mathcal{B}_k \cap \mathcal{B}_{k+1}^c$ implies 
that there is a hierarchy in $\mathcal{B}_k'$ but,
by \textit{3.}~of Definition \ref{Definition: hierarchy}, there is no 
edge between $\mathcal{D}_{N_k}^{(\rho,\ell)}(z)$ 
and $\mathcal{D}_{N_k}^{(\rho,\ell)}(z')$ for any pair $(z,z')$ 
as in the definition of $\mathcal{A}_2$. Hence, 
$\mathcal{B}_k \cap \mathcal{B}_{k+1}^c \subset \mathcal{B}_k' \cap \mathcal{A}_2$,  
and therefore we obtain  
\begin{equation}\label{1,2}
\p_{x,y}\left[\mathcal{B}_k \cap \mathcal{B}_{k+1}^c\right]
\le
\p_{x,y}\left[\mathcal{B}_k' \cap \mathcal{A}_2\right]
\le 
\p_{x,y}\left[\mathcal{B}_k' \cap \mathcal{A}_1\right] + 
\p_{x,y}\left[\mathcal{B}_k' \cap  \mathcal{A}_1^c \cap \mathcal{A}_2\right], 
\end{equation}
and it remains to bound the two terms on the right-hand side. 
To bound the first term we note that 
for $N$ sufficiently large it holds that for any hierarchy in $\mathcal{B}_k'$ and 
$\sigma, \sigma' \in \{0,1\}^k$ with $z_\sigma\neq z_{\sigma'}$, 
$R_{N_k+\ell}(z_\sigma)\cap R_{N_k+\ell}(z_{\sigma'})=\emptyset$ 
for all $\ell\in(N_m,2N_m)$. 
The latter implies that the events 
$\left\{|\mathcal{D}_{N_k}^{(\rho,\ell)}(z_\sigma)| \le  \rho_0 (2N_k)^d\right\}$ 
are independent for different $\sigma\in\{0,1\}^k$. 
It follows that for $N$ sufficiently large,
\begin{eqnarray*}
\p_{x,y}\left[\mathcal{B}_k' \cap \mathcal{A}_1\right] 
&\le& 
2^k
\p\left[|\mathcal{D}_{N_k}^{(\rho,\ell)}(0)| \le  \rho_0 (2N_k)^d\right].
\end{eqnarray*}
By Corollary \ref{Corollary:3.3&3.4} $(ii)$ 
there exist $\rho_0 \in (0,1)$ and $\ell_0 < \infty$ such that 
for all $\ell\in(\ell_0,N_k/\ell_0)$, 
\begin{equation*}
\p\left[
|\mathcal{D}_{N_k}^{(\rho_0,\ell)}(0)| \le  \rho_0 (2N_k)^d 
\right] 
\le 
e^{-\rho_0 N_k^{2d-\alpha'}} 
\end{equation*}
(although the definition of $\mathcal{D}_{N_k}^{(\rho_0,\ell)}(0)$ differs 
from $\mathcal{D}_{N_k}^{(\rho_0,\ell)}$ in 
Corollary \ref{Corollary:3.3&3.4} $(ii)$ we can still use the 
result because $R_{N_k}(0)$ contains a box 
of side length $N_k/3$). 
If we choose $N$ so large that 
$N_m \ge \ell_0$ and $N_{k}/\ell_0 \ge 2N_m$ for all $k=1, \ldots, m-1$, 
we finally obtain for all N sufficiently large, $\rho \in (0, \rho_0)$ 
and $\ell\in(N_m,2N_m)$,  
\begin{equation}\label{Equation: ell}
\p_{x,y}\left[\mathcal{B}_k' \cap \mathcal{A}_1\right] 
\le
2^k e^{-\rho_0 N_k^{2d-\alpha'}}
\le
2^k e^{-\rho_0 N_m^{2d-\alpha'}}
\le
2^k e^{-\rho_0 N_m^{2d\gamma-\alpha'}}.
\end{equation} 
We proceed as in the derivation of \eqref{Equation: trick}, 
using the same arguments as in \cite{biskup}, to see that 
the second term in \eqref{1,2} satisfies  
\begin{equation*}
\p_{x,y}\left[\mathcal{B}_k' \cap  \mathcal{A}_1^c \cap \mathcal{A}_2\right] 
\le
2^{k}
\exp\left\{-\lambda \rho_0^22^{2d} N_k^{2d} (5dN_{m-1})^{-\alpha} \right\}
\le
2^{k}
\exp\left\{-\lambda \rho_0^22^{2d}  (5d)^{-\alpha} N_{m-1}^{2d\gamma-\alpha} \right\}, 
\end{equation*}
for all $N$ sufficiently large. 
Using \eqref{1,2}, \eqref{Equation: ell} and that 
$N_{m-1}^{2d\gamma-\alpha} \ge N_{m}^{2d\gamma-\alpha'}$, 
it follows that for all $N$ sufficiently large, 
\begin{equation*}
\p_{x,y}\left[\mathcal{B}_k \cap \mathcal{B}_{k+1}^c\right]
\le 
2^k e^{-\rho_0 N_m^{2d\gamma-\alpha'}} +
2^{k}
\exp\left\{-\lambda \rho_0^22^{2d}  (5d)^{-\alpha} N_{m-1}^{2d\gamma-\alpha} \right\}
\le 
2^{k+1}e^{ -c_6N_m^{2d\gamma-\alpha'}}, 
\end{equation*}
with $c_6 = \min\{\rho_0 , \lambda \rho_0^22^{2d}  (5d)^{-\alpha} \}$. 

~

To prove \textit{(iii)} we note that 
any hierarchy in $\mathcal{B}_m$ has $2^{m-1}$ gaps 
whose $2^m$ particles $\{v_1=x,v_2=y,v_3, \ldots, v_{2^m}\}$ 
satisfy $\Lambda_\ell(v_i)\cap\Lambda_\ell(v_j)=\emptyset$ 
for all $i\neq j$.  
We therefore obtain for $N$ sufficiently large, 
\begin{eqnarray*}
\p_{x,y}\left[ \mathcal{B}_m \cap\mathcal{S}^c \right]
&\le&
2^m
\p_{x}\left[
X( \Lambda_\ell(x) )> e\nu (2\ell)^d \right].
\end{eqnarray*} 
Using Chernoff's bound, see \eqref{Equation: Chernoff}, and 
\eqref{Equation: (4.9)} we obtain 
for $\ell \ge N_m$, 
\begin{equation*}
2^m\p\left[X( \Lambda_\ell(0) )> e\nu (2\ell)^d \right]
~\le~
2^m
e^{-\nu(2\ell)^d}
~\le~
(\log N)^{\log 2/\log(1/\gamma)}
\exp\left\{-\nu2^de^{d(\log\log N)^\varepsilon}\right\},
\end{equation*}
which converges to $0$ as $N\to\infty$. 
\end{Proof}

~

\begin{Proof}[of Proposition \ref{Proposition: 4.1}]
Using part $(i)$ of Corollary \ref{Corollary:3.3&3.4} (to bound the 
probability of event $\mathcal{T}^c$)
and Lemma \ref{Lemma: 4.2,4.3}, 
Proposition \ref{Proposition: 4.1} 
is proven by using the same arguments as in \cite{biskup}.
\end{Proof}

\subsection{Finite variance of degree distribution, case 2}

In order to finish the proof of Theorem \ref{Theorem: graph distance} 
it remains to show statement \textit{(b2)}. 
We start with the following lemma.

\begin{lemma}\label{Lemma: edgesize}
Assume $\min\{\alpha, \beta\alpha\} > d$. 
For all $\delta \in (0,\alpha(\beta \wedge 1)-d)$ 
there exist $t_0 \ge 1$ and a constant $c_7>0$ such that for 
all $s\ge1$ and $t \ge t_0$,
\begin{equation*}
\p\left[\text{there is an edge in $\Lambda_s$ with size at least $t$} \right] 
\le 
c_7 s^d t^{d-\alpha( \beta \wedge 1)+\delta}.
\end{equation*}
\end{lemma}

~

\begin{Proof}[of Lemma \ref{Lemma: edgesize}]
Using \eqref{Equation: Pareto} we obtain for all $|x-y|\ge t_0$ 
with $t_0$ sufficiently large, 
\begin{eqnarray*}
\E\left[\left(\lambda\frac{W_xW_y}{|x-y|^\alpha}\right) \wedge 1\right]
&\le& 
(1+1_{\{\beta\neq1\}}/|\beta-1|)\lambda^{\beta \wedge 1}\left( 1+ \max\{1, \beta\} \log (\lambda^{-1}|x-y|^\alpha)\right )^2 |x-y|^{-\alpha(\beta \wedge 1)}
\\&\le&
|x-y|^{-\alpha(\beta \wedge 1)+\delta}.
\end{eqnarray*}
It follows that for all $t\ge t_0$, 
using $1-e^{-x}\le x\wedge1$, 
\begin{eqnarray*}
&&\hspace{-0.5cm}
\p\left[\text{there is an edge  in $\Lambda_s$ with size at least $t$} \right] 
~\le~
\E\left[
\sum_{x,y \in X \cap \Lambda_s}
1_{\{|x-y|>t\}}
|x-y|^{-\alpha(\beta \wedge 1)+\delta}
\right]
\\&&=
\sum_{k\ge1}\frac{\p\left[X\left(\Lambda_s\right)=k\right]}{(2s)^{dk}}
\int_{\Lambda_s}\cdots\int_{ \Lambda_s}
\sum_{i=1}^k\sum_{j=1}^k1_{\{|x_i-x_j|>t\}}|x_i-x_j|^{-\alpha(\beta \wedge 1)+\delta}
dx_1\cdots dx_k
\\&&=
\sum_{k\ge1}\frac{\p\left[X\left(\Lambda_s\right)=k\right]}{(2s)^{2d}}k(k-1)
\int_{\Lambda_s}\int_{ \Lambda_s}
1_{\{|x-y|>t\}}|x-y|^{-\alpha(\beta \wedge 1)+\delta}
dx\, dy
\\&&\le
\sum_{k\ge1}\frac{\p\left[X\left(\Lambda_s\right)=k\right]}{(2s)^{d}}k(k-1)
\int_{|y| > t}|y|^{-\alpha(\beta \wedge 1)+\delta}dy
~=~
\nu^2(2s)^d\int_{|y| > t}|y|^{-\alpha(\beta \wedge 1)+\delta}dy,
\end{eqnarray*}
where in the last step we used that $X(\Lambda_s)$ has a Poisson 
distribution with parameter $\nu (2s)^d$. 
It follows that for an appropriate constant $c_7>0$ and for all $t\ge t_0$ 
with $t_0$ sufficiently large, 
\begin{equation*}
\p\left[\text{there is an edge  in $\Lambda_s$ with size at least $t$} \right]
~\le~
c_7 s^d t^{d-\alpha( \beta \wedge 1)+\delta},
\end{equation*}
which finishes the proof of Lemma \ref{Lemma: edgesize}. 
\end{Proof}

~

\begin{Proof}[of \textit{(b2)} of Theorem \ref{Theorem: graph distance}]
Once the proof of Lemma \ref{Lemma: edgesize} is established, 
the proof of \textit{(b2)} of Theorem \ref{Theorem: graph distance} 
follows one-to-one from the proof of 
Theorem 1 of \cite{Berger2} and  Theorem 8 \textit{(b2)} of \cite{Rajat}. 
There, a renormalization is applied to see that 
we have a linear lower bound on the graph distances 
within ``good'' finite boxes, see also 
Definition 2 and Lemma 2 of \cite{Berger2}. 
Lemma \ref{Lemma: edgesize} is then used to prove 
that all centered boxes of sufficiently large side lengths 
are ``good'', a.s., see Lemma  14 of \cite{Rajat}, 
which then implies  \textit{(b2)} of Theorem \ref{Theorem: graph distance}. 
We refer to \cite{Rajat} for the details of the proof. 
\end{Proof}


\end{document}